\documentclass[11pt,a4paper]{amsart}
\usepackage[latin1]{inputenc}
\usepackage[T1]{fontenc}
\usepackage[english]{babel}

\usepackage{amsmath,amssymb,amsthm,amsfonts}

\usepackage{latexsym}

\usepackage{fancybox}
\usepackage{diagrams}

\theoremstyle{plain}
\newtheorem{theorem}{Theorem}[section]
\newtheorem{lemma}[theorem]{Lemma}
\newtheorem{corollary}[theorem]{Corollary}
\newtheorem{prop}[theorem]{Proposition}

\theoremstyle{definition}
\newtheorem{definition}[theorem]{Definition}
\newtheorem{example}[theorem]{Example}
\newtheorem{examples}[theorem]{Examples}

\newtheorem{question}[theorem]{Question}
\newtheorem{questions}[theorem]{Questions}

\newtheorem{remark}[theorem]{Remark}
\newtheorem{remarks}[theorem]{Remarks}

\newcommand{\C}{\mathbb{C}}

\newcommand{\R}{\mathbb{R}}
\newcommand{\K}{\mathbb{K}}
\newcommand{\N}{\mathbb{N}}
\newcommand{\nat}{\mathbb{N}}
\newcommand{\T}{\mathbb{T}}

\newcommand{\eps}{\varepsilon}
\newcommand{\ext}[1][X^*]{\ensuremath{\mathrm{ext}(B_{#1})}}
\newcommand{\extr}{\ensuremath{\mathrm{ext}}}
\newcommand{\Id}{\mathrm{Id}}
\newcommand{\diam}{\mathrm{diam}}
\newcommand{\dist}{\mathrm{dist}}

\newcommand{\Ker}{\mathrm{Ker}}

\newcommand{\conv}{\mathrm{conv}}
\newcommand{\cconv}{\overline{\mathrm{conv}}}
\newcommand{\aconv}{\mathrm{aconv}}
\newcommand{\caconv}{\overline{\mathrm{aconv}}}

 \DeclareMathOperator{\re}{Re}
 
 \DeclareMathOperator{\dent}{dent}

 \renewcommand{\leq}{\leqslant}
\renewcommand{\geq}{\geqslant}

\begin{document}

\begin{center}
[Trans. Amer. Math. Soc. (to appear)]
\end{center}

\title[Slicely Countably Determined Banach spaces]{Slicely Countably Determined Banach
spaces}
\author[Avil\'{e}s]{Antonio Avil\'{e}s}
\author[Kadets]{Vladimir Kadets}
\author[Mart\'{\i}n]{Miguel Mart\'{\i}n}
\author[Mer\'{\i}]{Javier Mer\'{\i}}
\author[Shepelska]{Varvara Shepelska}

\address[Avil\'{e}s]{Departamento de Matem\'{a}ticas\\ Universidad de
Murcia\\ 30100 Murcia Spain} \email{\texttt{avileslo@um.es}}
\address[Kadets \& Shepelska]{Department of Mechanics and Mathematics,
Kharkov National University, pl.~Svobody~4,  61077~Kharkov,
Ukraine} \email{\texttt{vova1kadets@yahoo.com} \qquad
\texttt{shepelskaya@yahoo.com}}
\address[Mart\'{\i}n \& Mer\'{\i}]{Departamento de An\'{a}lisis Matem\'{a}tico \\ Facultad de
 Ciencias \\ Universidad de Granada \\ 18071 Granada, Spain}
\email{\texttt{mmartins@ugr.es} \qquad \texttt{jmeri@ugr.es} }


\date{September 15th, 2008. Revised: February 9th, 2009}

\thanks{First named author supported by Marie Curie Intra-European
Fellowship MCEIF-CT2006-038768 (EU) and Spanish research project
MTM2005-08379 (MEC and FEDER). Second named author supported by
Junta de Andaluc\'{\i}a and FEDER grant P06-FQM-01438. Third and fourth
named authors partially supported by Spanish MEC and FEDER project
no.\ MTM2006-04837 and Junta de Andaluc\'{\i}a and FEDER grants FQM-185
and P06-FQM-01438. Fifth named author was partially supported by
N.~I.~Akhiezer foundation.}

\begin{abstract}
We introduce the class of slicely countably determined Banach
spaces which contains in particular all spaces with the RNP and
all spaces without copies of $\ell_1$. We present many examples
and several properties of this class. We give some applications to
Banach spaces with the Daugavet and the alternative Daugavet
properties, lush spaces and Banach spaces with numerical
index~$1$. In particular, we show that the dual of a real
infinite-dimensional Banach with the alternative Daugavet property
contains $\ell_1$ and that operators which do not fix copies of
$\ell_1$ on a space with the alternative Daugavet property satisfy
the alternative Daugavet equation.
\end{abstract}

\maketitle

\section{Introduction}
The aim of this paper is to introduce the class of slicely
countably determined Banach spaces, give many examples and several
properties of this class and, finally, to use this concept to give
some applications to Banach spaces with the Daugavet property and
to Banach spaces with numerical index~$1$. Let us introduce the
needed notation and definitions.

Given a Banach space over $\K$ ($\K=\R$ or $\K=\C$), we write
$S_X$ for its unit sphere and $B_X$ for its closed unit ball. The
dual space of $X$ is denoted by $X^*$ and $L(X)$ is the Banach
algebra of all bounded linear operators from $X$ to $X$. The space
$X$ has the \emph{Daugavet property} \cite{KadSSW} if every
rank-one operator $T\in L(X)$ satisfies
\begin{equation}\label{DE}\tag{\textrm{DE}}
\|\Id + T\|=1 + \|T\|.
\end{equation}
In this case, all operators on $X$ which do not fix copies of
$\ell_1$ (in particular, weakly compact operators) also satisfy
\eqref{DE} \cite{Shv1}. If every rank-one operator $T\in L(X)$
satisfies the norm equality
\begin{equation}\label{aDE}\tag{\textrm{aDE}}
\max_{\theta\in\T}\|\Id + \theta\,T\|=1 + \|T\|
\end{equation}
($\T$ being the set of modulus one scalars), $X$ has the
\emph{alternative Daugavet property} \cite{MaOi} and then all
weakly compact operators on $X$ also satisfy \eqref{aDE}. A Banach
space is said to have \emph{numerical index~$1$} \cite{D-Mc-P-W}
if every $T\in L(X)$ satisfies that $v(T)=\|T\|$, where
$$
v(T)=\bigl\{|x^*(Tx)|\ : \ x\in S_X,\,x^*\in S_{X^*},\, x^*(x)=1\bigr\}
$$
is the \emph{numerical radius} of the operator $T$. It is known
\cite{D-Mc-P-W} that
$$
v(T)=\|T\| \qquad \Longleftrightarrow \qquad T \text{ satisfies \eqref{aDE}.}
$$
Then, $X$ has numerical index~$1$ if and only if every $T\in L(X)$
satisfies \eqref{aDE}. It follows from the above discussion that
\begin{diagram}
\Ovalbox{Daugavet property} & \rImplies & \Ovalbox{Alternative
Daugavet property} & \lImplies & \Ovalbox{Numerical index~$1$}
\end{diagram}
None of the above implications reverses in general
\cite[Example~3.2]{MaOi}. For the first implication, it is even
known that it is not reversible under any isomorphic property
\cite[Corollary~3.3]{MaOi}. On the other hand, it is known that
the second implication reverses for Asplund spaces and for Banach
spaces with the Radon-Nikod\'{y}m property \cite[Remark~6]{LMP1999}.
We refer the interested reader to \cite{I-K-W,KMMP,KaMaPa,Mar-ADP}
and the already cited references for recent results, more
information and background on these properties.

We will say that $X$ is \emph{slicely countably determined}
(\emph{SCD} in short) if every bounded convex subset $A$ of $X$ is
an \emph{SCD set}, i.e.\ there is a sequence $\{S_n\,:\,n\in\N\}$
of slices of $A$ such that $A\subseteq \cconv(B)$ whenever
$B\subseteq A$ intersects all the $S_n$'s. Here a \emph{slice} of
a convex set $A$ is the subset given by
$$
S(A,x^*,\eps)=\{x\in A\,:\, \re x^*(x)>\sup\re x^*(A)-\eps\}
$$
and $\cconv(\cdot)$ stands for the closed convex hull. This
isomorphic property, which clearly implies separability, is
sufficient to get numerical index~$1$ from the alternative
Daugavet property and it is weaker than both RNP and being Asplund
(for separable spaces). Actually, this property is satisfied by
both separable strongly regular spaces and separable Banach spaces
which do not contain copies of $\ell_1$. This is the main
motivation of the study of SCD spaces.

In section~\ref{sec:sets} we study SCD sets, giving examples and
elementary properties. We show, for instance, that the sequence of
slices can be replaced by a sequence of relatively weakly open
sets or by a sequence of convex combinations of slices. In
section~\ref{sec:spaces} we study SCD spaces and show some
stability properties. For instance, it is a three space property,
so it is stable for finite sums, and it is stable for some
infinite unconditional sums.

Since it is not easy to deal with Banach spaces with numerical
index~$1$, there are in the literature several geometrical
sufficient conditions (see \cite{KaMaPa}), the weakest one being
the so-called lushness. A Banach space $X$ is said to be
\emph{lush} \cite{BKMW} if for every $x,y\in S_X$ and every
$\eps>0$, there is a slice $S=S(B_X,x^*,\eps)$ with $x^* \in
S_{X^*}$ such that $x \in S$ and $\dist \left(y,\aconv (S)\right)
< \eps$ (where $\aconv(A)$ denotes the absolutely convex hull of
the set $A$). Lush spaces have numerical index~$1$
\cite[Proposition~2.2]{BKMW}, but it has been very recently shown
that the converse result is not true \cite{KMMS}. We refer to
\cite{BKMM-lush,BKMW} for background.

It is actually shown in section~\ref{sec:applications-to-ni1} that
an SCD Banach space with the alternative Daugavet property is
lush. This result allows us to show that $\ell_1$ embeds in the
dual of every real infinite-dimensional Banach space with the
alternative Daugavet property. This answers in the positive
\cite[Problem~18]{KaMaPa}.

Section~\ref{sec:SCD-operators} is devoted to SCD-operators and
hereditary-SCD-operators. A bounded linear operator
$T:X\longrightarrow Y$ between two Banach spaces $X$ and $Y$ is
said to be an \emph{SCD-operator} if $T(B_X)$ is an SCD set, and
$T$ is a \emph{hereditary-SCD-operator} if every bounded convex
subset of $T(B_X)$ is SCD. We show that SCD-operators on a Banach
space with the alternative Daugavet property satisfy \eqref{aDE}.
Therefore, operators which do not fix copies of $\ell_1$ on a
Banach space with the alternative Daugavet property satisfy
\eqref{aDE}. For a Banach space with the Daugavet property it is
shown that every SCD-operator is strong Daugavet (and so it
satisfies \eqref{DE}), and every hereditary-SCD-operator is
narrow.

Section~\ref{sec:countable-pi-base} is devoted to the study of
sets with a countable $\pi$-base of the weak topology. It is shown
in section~\ref{sec:sets} that these sets are SCD, but it is not
known whether the converse result is true. It is also shown in
section~\ref{sec:sets} that separable sets without $\ell_1$
sequences have countable $\pi$-bases of the weak topology, and in
this section we show that the same is true for CPCP sets and for
bounded convex subsets of both $c_0(\ell_1)$ and $\ell_1(c_0)$. We
also show some characterizations of SCD sets which remind of the
existence of countable $\pi$-bases of the weak topology. One of
these characterizations allows us to show that the set of extreme
points of the weak$^*$-closure (in the bidual space) of an SCD set
has a countable $\pi$-base of the weak$^*$ topology, and so it is
weak$^*$ separable. The set of extreme points of a convex set $B$
will be  denoted by $\extr(B)$.

Finally, section~\ref{sec:openquestions} contains several open
questions.

\section{Slicely countably determined sets}\label{sec:sets}

\begin{definition}\label{SCD sets-def 1}
Let $X$ be a Banach space and let $A$ be a convex bounded subset
of $X$. A countable family $\{V_n \,:\, n\in\nat\}$ of subsets of
$A$ is called \emph{determining} for $A$ if $A \subseteq \cconv
(B)$ for every $B\subseteq A$ intersecting all the sets $V_n$.
Equivalently, $\{V_n \,:\, n\in\nat\}$ is determining for $A$ if
for every sequence $\{v_n\}_{n\in\nat}$ with $v_n\in V_n$
($n\in\nat$), one has $A \subseteq
\cconv\bigl(\{v_n\,:\,n\in\nat\}\bigr).$
\end{definition}

We give three easy observations which will be useful later on. The
first one is a consequence of the Hahn-Banach theorem. The second
and third ones are straightforward.

\begin{prop}\label{determ-sequences}
Let $X$ be a Banach space and let $A$ be a convex bounded subset
of $X$. A sequence $\{V_n\,:\,n\in\N\}$ of subsets of $A$ is
determining if and only if every slice of $A$ contains one of the
$V_n$.
\end{prop}

\begin{proof}
The ``if'' part is evident: if $B\subseteq A$ intersects all the
$V_n$, then it intersects all the slices of $A$, and then by the
Hahn-Banach theorem $\cconv(B)\supseteq A$. Now the ``only if''
part. Assume that some slice $S$ of $A$ does not contain any of
the $V_n$. Then $A \setminus S$ is a convex relatively closed
subset of $A$ intersecting all the $V_n$. But $A \setminus S \neq
A$, which means that $\{V_n\,:\,n\in\N\}$ is not determining.
\end{proof}

\begin{remark}
{\slshape Let $X$ be a Banach space and let $A$ be a convex
bounded subset of $X$. Suppose that there is a sequence $\{a_n\, :
\, n\in \N\}$ of points in $A$ such that $A\subseteq
\cconv\bigl(\{a_n\, : \, n\in \N\}\bigr)$ and that for every
$n\in\N$, there is sequence $\{V_{n,m}\,:\, m\in\N\}$ of subsets
of $A$ such that $a_n\in \cconv(B)$ whenever $B\subseteq A$
intersects $V_{n,m}$ for every $m\in\N$. Then, the family
$\{V_{n,m}\,:\, n,m\in\N\}$ is determining for $A$.\ }
\end{remark}

As an immediate consequence of the above result, we get the
following.

\begin{remark}\label{rem:dense-sequence}
{\slshape Let $X$ be a Banach space and let $A$ be a separable
convex bounded subset of $X$. Suppose that for every $a\in A$
there is a sequence $\{V_{m}^{a}\,:\, m\in\N\}$ of subsets of $A$
such that $a\in \cconv(B)$ whenever $B\subseteq A$ intersects
$V_{m}^{a}$ for every $m\in\N$. Then, taking a dense sequence
$\{a_n\,:\,n\in\N\}$ in $A$, the family $\{V_{m}^{a_n}\,:\,
n,m\in\N\}$ is determining for $A$.\ }
\end{remark}

We can now give the main definition of this section.

\begin{definition}\label{def:SCDset}
A convex bounded subset $A$ of a Banach space $X$ is said to be
\emph{slicely countably determined} (\emph{SCD set} in short) if
there is a determining sequence of slices of $A$.
\end{definition}

Two remarks are pertinent.

\begin{remark}
{\slshape It is clear from the definition that every SCD set is
separable.\ }
\end{remark}

\begin{remark}\label{SCD sets-rem 3}
{\slshape A convex bounded subset $A$ of a Banach space $X$ is SCD
if and only if the closure of $A$ is an SCD set.\ }
\end{remark}

\begin{proof}
Let us show first that $\overline A$ is SCD when $A$ is. Consider
a determining sequence of slices $S_n= S(A,x_n^*, \eps_n)$
($n\in\N$) for $A$, and let us prove that the slices $S'_n=
S({\overline A},x_n^*, \eps_n/2)$ ($n\in\N$) form a determining
sequence for the closure of $A$. Consider an arbitrary slice $S=
S({\overline A},x^*, \eps)$ of $\overline A$. Then, $S({\overline
A},x^*, \eps/2) \cap A = S(A, x^*, \eps/2)$ is a slice of $A$, so
there is $n\in\N$ such that $S(A,x^*, \eps/2) \supseteq S_n$ by
Proposition~\ref{determ-sequences}. Therefore, $S$ contains the
closure of $S_n$, which in turn contains $S'_n$, and again
Proposition~\ref{determ-sequences} gives us that $\{S'_n\}$ is
determining for ${\overline A}$.

For the converse implication, we consider a determining sequence
$\{S({\overline A},x_n^*, \eps_n)\,:\,n\in\N\}$ for $\overline A$,
and it is straightforward to show that $\{S(A,x_n^*,
\eps_n)\,:\,n\in\N\}$ is determining for $A$.
\end{proof}

Our fist goal is to present the basic examples related to
Definition~\ref{def:SCDset}: Radon-Nikod\'{y}m and Asplund sets are
SCD, whereas the unit ball of a Banach space with the Daugavet
property is not.

We start with subsets having sufficiently many denting points. Let
$X$ be a Banach space and let $A$ be a closed convex bounded
subset of $X$. A point of $A$ is said to be a \emph{denting point}
if it belongs to slices of $A$ of arbitrarily small diameter. We
write $\dent(A)$ to denote the set of denting points of $A$. We
say that $A$ is \emph{dentable} (in the sense of
Ghoussoub-Godefroy-Maurey-Schachermayer \cite[\S
III]{Gho-God-Mau-Scha}) if $A=\cconv\bigl(\dent(A)\bigr)$
\cite[Proposition~III.3]{Gho-God-Mau-Scha}.

\begin{prop}\label{prop:denting}
Let $X$ be a Banach space and let $A$ be a closed convex bounded
subset of $X$. If $A$ is separable and dentable, then $A$ is SCD.
\end{prop}

\begin{proof}
Since $A$ separable, so is the set of its denting points, so we
may find a countable collection of denting points $\{a_n\,:\,
n\in\nat\}$ of $A$ which is dense in $\dent(A)$. Now, for every
$n,m\in \nat$, we consider a slice $S_{n,m}$ of $A$ containing
$a_n$ and having diameter less than $1/m$. Then, the sequence
$\{S_{n,m}\,:\, n,m\in\nat\}$ is determining for $A$. Indeed, if
$B\subseteq A$ intersects all the $S_{n,m}$, then $a_n\in
\overline{B}$ for every $n\in \nat$, so
\begin{equation*}
A\subseteq \cconv\bigl(\dent(A)\bigr)=\cconv\bigl(\{a_n\,:\,n\in\nat\}\bigr)
\subseteq \cconv(\overline{B})=\cconv(B).\qedhere
\end{equation*}
\end{proof}

We recall that there is a concept of \emph{Radon-Nikod\'{y}m set}
(defined in terms of vector measures) which is equivalent to
dentability of all its closed convex bounded subsets (see \cite[\S
5]{Ben-Lind} or \cite[\S 2]{Bourgin}).

\begin{example}\label{example:RNP-sets}
{\slshape Let $X$ be a Banach space and let $A$ be a closed convex
bounded separable Radon-Nikod\'{y}m subset of $X$. Then, $A$ is an SCD
set.\ }
\end{example}

The norm $\|\cdot\|$ on a Banach space $X$ is said to be
\emph{LUR} at $x_0\in S_X$, if $\lim \|x_n - x_0\|=0$ whenever
$(x_n)_{n\in\N}\subseteq B_X$ is such that $\lim \|x_n+x_0\|=2$.
If the norm is LUR at each point of $S_X$, we say that $X$ (or its
norm) is \emph{LUR} (see \cite[Chapter~II]{DGZ} for background).
It is clear that every point in the unit sphere of a Banach space
$X$ with a LUR norm is denting so, in this case, $B_X$ is
dentable.

\begin{example}\label{example:LUR=>SCD}
{\slshape Let $X$ be a separable Banach space with a LUR norm.
Then, $B_X$ is SCD.\ }
\end{example}

It is well known that every separable Banach space admits a LUR
renorming (see \cite[Theorem~II.2.6.]{DGZ}). Therefore, the
following result follows immediately from
Proposition~\ref{prop:denting}.

\begin{example}\label{example:renorming-SCD}
{\slshape Every separable Banach space $X$ admits an equivalent
norm $|\cdot|$ such that $B_{(X,|\cdot|)}$ is an SCD set.\ }
\end{example}

Our second family of elementary examples of SCD sets deals with
the so-called \emph{Asplund property}, a concept related to
differentiability of convex continuous functions, which can be
equivalently reformulated in terms of separability and duality
\cite[\S5]{Bourgin}. A separable closed convex bounded subset $A$
of a Banach space $X$ has the Asplund property if and only if the
semi-normed space $(X^*,\rho_A)$ is separable, where
$$
\rho_A(x^*)=\sup\{|x^*(a)|\,:\, a\in A\} \qquad (x^*\in X^*).
$$
Of course, separable closed convex bounded subsets of Asplund
spaces have the Asplund property.

\begin{example}\label{example:Asplund-sets}
{\slshape Let $X$ be a Banach space and let $A$ be a closed convex
bounded subset of $X$. If $A$ is separable and has the Asplund
property, then, $A$ is SCD.\ }
\end{example}

\begin{proof}
We take a $\rho_A$-dense countable family $\{x_n^*\,:\, n\in\N\}$
in $(X^*,\rho_A)$, and consider the slices
$$
S_{n,m}=S(A,x_n^*,1/m) \qquad (n,m\in\N).
$$
We are done by just proving that if $\{v_{n,m}\,:\, n,m\in \N\}$
satisfies that $v_{n,m}\in S_{n,m}$ for every $n,m\in\N$, then
$$
A\subseteq \cconv \left(\{v_{n,m}\,:\,n,m\in\N\}\right).
$$
Indeed, suppose to the contrary that there are $a\in A$, $x^*\in
X^*$, and $\delta>0$ such that
$$
\re x^*(a) > \sup_{n,m} \re x^*(v_{n,m}) + \delta.
$$
Now, we may find $N\in\N$ such that $\rho_A(x_N^*-x^*)<\delta/2$
and so
\begin{multline*}
\re x_N^*(a) + \delta/2 > \re x^*(a) >\sup_{n,m}\re x^*(v_{n,m}) +
\delta
\\ \geq \sup_{m} \re x^*(v_{N,m}) + \delta > \sup_m \re x_N^*(v_{N,m}) + \delta/2
 = \sup \re x_N^*(A) + \delta/2,
\end{multline*}
a contradiction.
\end{proof}

We now show that there are convex bounded subsets of separable
Banach spaces which are not SCD.

\begin{example}\label{Examples-Daugavet}
{\slshape Let $X$ be a separable Banach space with the Daugavet
property. Then, $B_X$ is not an SCD set. In particular,
$B_{C[0,1]}$ and $B_{L_1[0,1]}$ are not SCD sets.\ }
\end{example}

\begin{proof}
Fix $x_0\in S_X$ and an arbitrary sequence of slices $(S_n)_{n \in
\N}$. We will get the result by showing that there is a sequence
$(x_n)_{n\in\N}$ such that $x_n\in S_n$ for every $n\in\N$ and
such that $x_0\notin \overline{\textrm{lin}\{x_n\,:\,n \in \N\}}$.
To do so, we use \cite[Lemma~2.8]{KadSSW} which says, in
particular, that for every finite-dimensional subspace $Y
\subseteq X$, every $\eps > 0$, and every slice $S$ of $B_{X}$,
there is an $x \in S$ such that
\begin{equation*}
\|y+tx\|\geq (1-\eps)(\|y\|+|\,t|) \qquad\forall y\in Y.
\end{equation*}
Using this result, one can select inductively elements $x_n \in
S_n$, $n\in \N$, in such a way, that
$$
\|y+tx_n\|\geq \left(1-\frac{1}{4^n}\right)(\|y\|+|\,t|) \qquad \bigl(y\in
\textrm{lin}\{x_k\,:\, k < n\}\bigr).
$$
Then, $\{x_n\,:\,n=0,1,\ldots\}$ form a sequence equivalent to the
unit vector basis of $\ell_1$, so $x_0$ is not in the closure of
$\textrm{lin}\{x_n\,:\, n \in \N\}$, as desired.
\end{proof}

For the case of $C[0,1]$, it is possible to give a direct proof
without using the Daugavet property, which we include here for the
sake of completeness.

\begin{example}
{\slshape If $K$ is an uncountable metrizable compact space, then
the unit ball of $C(K)$ is not an SCD set. }
\end{example}

\begin{proof}
Let $\mathcal{M}$ be a maximal family of mutually orthogonal
continuous measures in $C(K)^*$. This induces a decomposition of
$C(K)^*$ as
\begin{equation}\label{eq-C01-noSCD-Antonio-1}
C(K)^* = \left[\bigoplus_{\mu\in \mathcal{M}}
L^1(\mu)\right]_{\ell_1}\oplus_1 \ell_1(K),
\end{equation}
where $\ell_1(K)$ is the family of all discrete measures (see
\cite[pp.~84--85]{Albiac-Kalton}, for instance). As a consequence,
we have
\begin{equation}\label{eq-C01-noSCD-Antonio-2}
C(K)^{**} =\left[\bigoplus_{\mu\in \mathcal{M}}
L^\infty(\mu) \right]_{\ell_\infty}\oplus_\infty \ell_\infty(K).
\end{equation}
Let us write the slices of $B_{C[0,1]}$ in the form
$$
U[\nu,\alpha]=\bigl\{x\in B_{C(K)}\, :\, \re\nu(x)>\alpha\bigr\},
$$
where $\nu\in C(K)^*$, $\|\nu\|=1$, and $-1 < \alpha <1$. Suppose,
for the sake of contradiction, that there existed a countable
family of slices $\mathcal{B}_0$ such that every other slice
contains one from the family. Then, for every $\mu\in
\mathcal{M}$, there exist
$V_\mu=U[\nu_\mu,\alpha_\mu]\in\mathcal{B}_0$ such that
$V_\mu\subseteq U[\mu,0]$.

Now, for each $\nu\in C(K)^*$ we write
$$
\mathrm{Supp}_\mathcal{M}(\nu) = \bigl\{\mu\in\mathcal{M}\, :\, \mu\not\perp
\nu\bigr\}.
$$
Notice that this is a countable set which corresponds to the
support of $\nu$ in the left-hand side of the decomposition
\eqref{eq-C01-noSCD-Antonio-1}. We claim that $\mu\in
\mathrm{Supp}_\mathcal{M}(\nu_\mu)$ for every $\mu\in\mathcal{M}$.
This leads to a contradiction with the facts that $\mathcal{B}_0$
and all the sets $\mathrm{Supp}_\mathcal{M}(\nu_\mu)$ are
countable, while $\mathcal{M}$ is uncountable. Let us prove the
claim. Suppose that $\mu\not\in
\mathrm{Supp}_\mathcal{M}(\nu_\mu)$ and let $g$ be an element of
the unit ball of $C(K)^{**}$ where $\nu_\mu$ attains its norm.
Consider $f\in L^{\infty}(\mu)$ the $\mu$-coordinate of $g$ when
we view $g$ as an element of the $\ell_\infty$-sum according to
\eqref{eq-C01-noSCD-Antonio-2}. Let now $g'$ be the element of
$C(K)^{**}$ obtained from $g$ by changing the $\mu$-coordinate
from $f$ to $-f$. This is a new element of the unit ball of
$C(K)^{**}$ which satisfies that $g'(\mu)=-g(\mu)$ while
$g'(\nu_\mu) = g(\nu_\mu) = 1$. Hence, for either $h=g$ or $h=g'$,
we have an element $h$ in the unit ball of $C(K)^{**}$ such that
$h(\nu_\mu)=1$ and $h(\mu)<0$. Since the unit ball of $C(K)$ is
dense in the unit ball of $C(K)^{**}$, it follows that
$V_\mu\setminus U[\mu,0]\neq\emptyset$.
\end{proof}

\begin{remark}
{\slshape A subset of an SCD set is not necessarily SCD.\ }
Indeed, let $X=C[0,1]$. By Example~\ref{example:renorming-SCD},
there is an equivalent norm $|\cdot|$ on $X$ such that
$A=B_{(X,|\cdot|)}$ is SCD. Now, it is possible to find
$\lambda>0$ such that $C=\lambda\,B_{(X,\|\cdot\|_\infty)}$ is
contained in $A$. Finally, $C$ is not SCD by
Example~\ref{Examples-Daugavet}.
\end{remark}

Our next goal is to extend the above preliminary examples to more
intriguing ones. We will use several times the so-called
\emph{Bourgain's lemma} \cite[Lemma~5.3]{Bour-Paris} (it was
rediscovered in \cite{Shv1}), so we state it for the sake of
completeness. We refer the reader to \cite[Lemma~7.3]{van-Duslt}
for a reference easier to get. We recall that a \emph{convex
combination of slices} of a convex bounded subset $A$ of a Banach
space $X$ is a subset of $A$ of the form $\sum\limits_{k=1}^m
\lambda_i\,S_i$ where $\lambda_i>0$, $\sum \lambda_i =1$ and the
$S_i$'s are slices of $A$.

\begin{lemma}[Bourgain's lemma]\label{lemma:Bourgain}
Let $X$ be a Hausdorff locally convex space and let $K\subseteq X$
be closed bounded and convex. Then, every nonempty relatively
weakly open subset of $K$ contains a convex combination of slices.
\end{lemma}

\begin{remark}
{\slshape The condition of closedness of the set in Bourgain's
lemma can be omitted.\ } Indeed, let $A$ be a convex bounded set
and let $U$ be a relatively weakly open subset of $A$. We denote
by $V$ a relatively weakly open subset of $\overline A$ such that
$V \cap A = U$. By Bourgain's lemma, there are slices $S_1, S_2,
\ldots, S_n$ of $\overline A$ and coefficients $\lambda_k
> 0$ of a convex combination, such that $\sum_1^n \lambda_k S_k
\subseteq V$. Then, $S_k \cap A$ are slices of $A$ and $\sum_1^n
\lambda_k S_k \cap A \subseteq V\cap A=U$.
\end{remark}

The first consequence is an easy observation.

\begin{prop}\label{SCD-sets-rem-2}
In the definition of SCD sets, instead of slices one can take
convex combinations of slices. Hence, by Bourgain's lemma above,
one can also take relatively weakly open subsets.
\end{prop}

\begin{proof}
Let $\{V_n\,:\,n\in\N\}$ be a determining sequence formed by
convex combination of slices of $A$. Now, for every $n\in\N$,
there exists a collection of slices $\{S_{n,m}\,:\,m=1,\ldots
k_n\}$ and positive numbers $\{\lambda_{n,m}\,:\,m=1,\ldots k_n\}$
with $\sum\limits_{m=1}^{k_n}\lambda_{n,m}=1$, such that
$\sum\limits_{m=1}^{k_n}\lambda_{n,m}S_{n,m}\subseteq V_n$. Then,
the collection of slices $\{S_{n,m}\,:\, n\in\N,\ 1\leq m\leq
k_n\}$ is determining for $A$. Indeed, let $B$ be a subset of $A$
such that $B\cap S_{n,m}\neq\emptyset$ for all $n,m$, and consider
$b_{n,m}\in B\cap S_{n,m}$ for every $n,m$. If we take
$a_n=\sum\limits_{m=1}^{k_n}\lambda_{n,m}b_{n,m}$, it is clear
that $a_n\in \conv(B)\cap V_n$. So we know that $\conv(B)\cap
V_n\neq\emptyset$ for all $n$, which by the assumption gives us
that $\cconv(B)\supseteq A$.

Finally, if $A$ has a determining sequence of relatively weakly
open subsets $\{V_n\,:\,n\in\N\}$, Bourgain's lemma allows us to
find convex combinations of slices inside the $V_n$'s and the
proof above shows that $A$ is SCD.
\end{proof}

The first consequence of this result is that
Proposition~\ref{prop:denting} can be extended from dentable sets
to \emph{huskable} sets (the same definition with relatively
weakly open sets instead of slices). With not much work, we are
going to extend the result to the following more general setting.
A closed convex bounded subset $A$ of a Banach space $X$ has
\emph{small combinations of slices}
\cite{Gho-God-Mau-Scha,Ros-1988} if every slice of $A$ contains
convex combinations of slices of $A$ with arbitrarily small
diameter.

\begin{theorem}
Let $X$ be a Banach space and let $A$ be a separable closed convex
bounded subset of $X$ having small combinations of slices. Then,
$A$ is an SCD set.
\end{theorem}

\begin{proof}
By \cite[Corollary~III.7]{Gho-God-Mau-Scha}, for every $x\in A$
and every $\eps>0$, there is a convex combination of slices of $A$
contained in $B(x,\eps)$. Now, we take a countable dense subset
$\{x_n\,:\,n\in\N\}$ of $A$ and for $(n,m)\in\N\times\N$, we take
$V_{n,m}$ a convex combination of slices of $A$ contained in
$B(x_n,1/m)$. Then, if $B\subseteq A$ intersects all the
$V_{n,m}$, it intersects also all the balls $B(x_n,1/m)$.
Therefore, the set $\{x_n\,:\,n\in\N\}$ is contained in $\overline
B$ and so, $A=\cconv(B)$. Finally,
Proposition~\ref{SCD-sets-rem-2} gives us that $A$ is SCD.
\end{proof}

RNP sets have small combinations of slices, so the above result
extends Example~\ref{example:RNP-sets}. Even more, strongly
regular sets (in particular, huskable sets, CPCP sets) have small
combinations of slices \cite[Proposition~III.5]{Gho-God-Mau-Scha}.
We recall that a closed convex bounded subset $A$ of a Banach
space is said to be \emph{strongly regular} if every non-empty
convex subset $L$ of $A$ contains a convex combination of slices
of $L$ of arbitrarily small diameter. $A$ has the \emph{convex
point of continuity property} (\emph{CPCP} in short) if every
convex closed subset $B$ of $A$ contains a weak-to-norm point of
continuity of the identity mapping. In this case, for every convex
subset $B$ of $A$ and for every $\eps>0$, there is a relatively
weakly open subset $C\subseteq B$ with $\diam(C)<\eps$
\cite{Bourgain1980}.

\begin{corollary}\label{SCDsets:corollary-strongly-regular}
Let $X$ be a Banach space and let $A$ be a closed convex bounded
subset of $X$. If $A$ is separable and strongly regular, then $A$
is SCD. In particular, separable CPCP sets are SCD.
\end{corollary}

Our next aim is to extend Example~\ref{example:Asplund-sets} to
sets which do not contain $\ell_1$ sequences. We need the
following topological definition. By a \emph{$\pi$-base} of a
topological space $(T,\tau)$ we understand a family $\{O_i\,:\,
i\in I\}$ of nonempty open sets such that every nonempty open
subset $O$ of $T$ contains one of the elements of the family. The
following result is another consequence of Bourgain's lemma.

\begin{prop} \label{SCD:prop-countablepibase}
Let $X$ be a Banach space and let $A$ be a convex bounded subset
of $X$. If $(A, \sigma(X,X^*))$ has a countable $\pi$-base, then
$A$ is an SCD set.
\end{prop}

\begin{proof}
Let $\{V_n\,:\, n\in \N\}$ be a countable $\pi$-base of $(A,
\sigma(X,X^*))$. Since slices of $A$ have non-empty weak interior,
any of them contains some of the $V_n$. But then,
Proposition~\ref{determ-sequences} shows that the sequence
$\{V_n\}$ is determining for $A$ and
Proposition~\ref{SCD-sets-rem-2} gives that $A$ is SCD.
\end{proof}

The main consequence of the above proposition is the following. We
recall that an \emph{$\ell_1$-sequence} of a Banach space is just
a bounded sequence which is equivalent to the natural basis of
$\ell_1$

\begin{theorem}\label{Examples-Theorem:ell1}
Let $X$ be a Banach space and let $A$ be a separable convex
bounded subset of $X$ which contains no $\ell_1$-sequences. Then,
$(A, \sigma(X,X^*))$ has a countable $\pi$-base. In particular,
$A$ is an SCD set.
\end{theorem}

\begin{proof}
By \cite[Theorem~3.11]{van-Duslt}, $(A,\sigma(X,X^*))$ is a
relatively compact subset of the space of first Baire class
functions on $(B_{X^*},\sigma(X^*,X))$, and we can apply
\cite[Lemma~4]{T} by Todor\v cevi\'c, to deduce that $(A,
\sigma(X,X^*))$ has a $\sigma$-disjoint $\pi$-base (i.e. a
$\pi$-base $\{V_i\, :\, i\in I\}$ such that
$I=\bigcup_{n\in\mathbb{N}}I_n$ and each subfamily $\{V_i\, :\,
i\in I_n\}$ is a pairwise disjoint family). Now, it is clear that
a $\sigma$-disjoint family of open subsets in a separable space
has to be countable. Finally, $A$ is SCD by
Proposition~\ref{SCD:prop-countablepibase}.
\end{proof}

This result obviously extends Example~\ref{example:Asplund-sets}
since Asplund sets cannot contain $\ell_1$-sequences.

\section{Slicely Countably Determined spaces}\label{sec:spaces}

\begin{definition}\label{SCD spaces-def 1}
A separable Banach space $X$ is said to be \emph{slicely countably
determined} (\emph{SCD space} in short) if every convex bounded
subset of $X$ is an SCD set.
\end{definition}

By just using the results of the previous section on SCD sets, we
get the main examples of SCD spaces.

\begin{examples}$ $\label{examples:main-SCD-spaces}
\begin{enumerate}
\item[(a)] {\slshape If $X$ is a separable strongly regular
    space, then $X$ is SCD. In particular, RNP spaces (more
    generally, CPCP spaces) are SCD.\ }
\item[(b)] {\slshape Separable spaces which do not contain
    copies of $\ell_1$ are SCD. In particular, if $X^*$ is
    separable, then $X$ is SCD.\ }
\item[(c)] {\slshape Both families include reflexive separable
    spaces, which are then SCD spaces.\ }
\end{enumerate}
\end{examples}

With respect to spaces which are not SCD, we only know of the
Daugavet spaces.

\begin{examples}$ $
\begin{enumerate}
\item[(a)] {\slshape If $X$ is a separable Banach space which
    admits an equivalent renorming with the Daugavet property,
    then $X$ is not SCD.\ }
\item[(b)] {\slshape In particular, there is a Banach space
    with the Schur property which is not an SCD space.\ }
    Indeed, in \cite{KadWer} the existence of a separable
    space having the Schur property and the Daugavet property
    at the same time was proved.
\end{enumerate}
\end{examples}

Let us state the following immediate observations.

\begin{remarks}$ $
\begin{enumerate}
\item[(a)] Every subspace of an SCD space is SCD.
\item[(b)] For quotients the situation is different. For
    instance, $C[0,1]$ is a non-SCD quotient of the SCD space
    $\ell_1$.
\end{enumerate}
\end{remarks}

Our next aim is to show some stability results for the SCD spaces.
The first one is a ``three space property''. We need the following
technical lemma which shows that in Definition~\ref{SCD spaces-def
1} it suffices to consider sets with nonempty interior.

\begin{lemma} \label{SCD-spaces-Lemma-open}
Let $X$ be a separable Banach space. If every \emph{open} convex
bounded subset of $X$ is SCD, then $X$ is SCD.
\end{lemma}

\begin{proof}
Our first observation is that our hypothesis forces that every
bounded convex subset $A$ of $X$ with nonempty interior is SCD.
Indeed, notice that since $A$ is convex, the closure of the
interior of $A$ coincides with the closure of $A$, and we may
apply Remark~\ref{SCD sets-rem 3} two times to get that $A$ is
SCD.

Now, let $A\subseteq X$ be bounded and convex. Since $X$ is
separable, we may find a sequence $\{x_n\,:\,n\in\N\}\subseteq A$
which is dense in $A$. Let $\{\varepsilon_n\}_{n\in\N}$ be a
sequence of positive reals which tends to zero. For every
$n,m\in\mathbb{N}$ fixed, we denote
$A_{n,m}=\conv\bigl(B_{\varepsilon_m}(x_n)\cup A\bigr)$ which
clearly contains $A$. Since the interior of $A_{n,m}$ is not
empty, we may find a determining sequence $\{S_{n,m}^k\,:\,k\in
\N\}$ of slices of $A_{n,m}$. Now, from the structure of
$A_{n,m}$, it follows that either $S_{n,m}^k\cap
B_{\varepsilon_m}(x_n)\neq\emptyset$, or $S_{n,m}^k\cap A
\neq\emptyset$. Let $K_{n,m}$ be the set of all indices $k\in\N$
for which $S_{n,m}^{k}$ intersects $A$, and denote
$\widetilde{S}_{n,m}^{k}=S_{n,m}^{k}\cap A$ for all $k\in
K_{n,m}$, which are clearly slices of $A$. Also note that for
every integer $k\notin K_{n,m}$, the slice $S_{n,m}^{k}$
intersects $B_{\varepsilon_m}(x_n)$. Finally, the family
$$
\left\{\widetilde{S}_{n,m}^{k}\ :\ n,m\in\N,\ k\in K_{n,m} \right\}
$$
is determining for $A$. Indeed, let $B$ be a subset of $A$
intersecting all the $\widetilde{S}_{n,m}^{k}$ and fix some
$\varepsilon>0$. Since the sequence $\{x_n\,:\,n\in\N\}$ is dense
in $A$, there is an integer $n_0\in\N$ and $b\in B$ such that
$\|b-x_{n_0}\|\leq \frac{\varepsilon}{2}$. Also, there is
$m_0\in\N$ such that $\varepsilon_{m_0}\leq
\frac{\varepsilon}{2}$, as $\varepsilon_m\rightarrow 0$ when
$m\rightarrow\infty$. We know that $B$ intersects all
$S_{n_0,m_0}^k$ with $k\in K_{n,m}$. On the other hand, we also
know that the slice $S_{n_0,m_0}^k$ intersects the ball
$B_{\varepsilon_{m_0}}(x_{n_0})$ for every $k\notin K_{n,m}$.
Hence we can deduce that the set $B_{n_0,m_0}=B\cup
B_{\varepsilon_{m_0}}(x_{n_0})\subseteq A_{n,m}$ intersects all
the $S_{n_0,m_0}^k$ which implies that
$$
\cconv \bigl(B_{n_0,m_0}\bigr)\supseteq A_{n_0,m_0}\supseteq A.
$$
Finally, notice that $B_{\varepsilon_{m_0}}(x_{n_0})\subseteq
B_{\frac{\varepsilon}{2}}(x_{n_0})\subseteq B_{\varepsilon}(b)$,
which implies that $B_{n_0,m_0}\subseteq B+\varepsilon B_{X}$.
Therefore, we can state that $\cconv \bigl(B+\varepsilon
B_{X}\bigr)\supseteq A$, and the arbitrariness of $\varepsilon$
gives us that $\cconv(B)\supseteq A$.
\end{proof}

We may now state the promised stability result.

\begin{theorem} \label{SCD spaces-thm2}
Let $X$ be a Banach space with a subspace $Z$ such that $Z$ and
$Y=X/Z$ are SCD spaces. Then, $X$ is also an SCD space.
\end{theorem}

\begin{proof}
We denote $q:X\longrightarrow Y=X/Z$ the quotient map. Let us show
that every \emph{open} convex bounded subset $A\subseteq X$ is
SCD, and then Lemma~\ref{SCD-spaces-Lemma-open} will imply that
$X$ is SCD. To do so, as $X$ is separable since $Y$ and $Z$ are,
and separability is a three-space property (see
\cite[Theorem~2.4.h]{CasGon}), we only need to find, for every
point $a\in A$, a sequence of weakly open subsets such that
whenever $B\subseteq A$ intersects every member of the sequence,
then $a\in \cconv(B)$ (see Remark~\ref{rem:dense-sequence}). We
fix some $a\in A$ and denote $A_a=\{x\in A\, :\, q(x)=q(a)\}$.
Then, $A_a$ is affine isomorphic to an open convex bounded subset
of $Z$ which is an SCD space (indeed, $A_a=(Z+a) \cap A$). It
follows that there is a determining sequence $\{S_n\}$ of slices
of $A_a$. Let $\{\widetilde{S}_n\}$ be their extensions to $A$.
For every $n\in\N$, consider $q(\widetilde{S}_n)\subseteq Y$,
which is open bounded and convex (its openness is a consequence of
the Open Mapping Theorem). Now, as long as $Y$ is SCD, we may find
a determining sequence $\{S_{n,m}\,:\, m\in\N\}$ of slices of
$q(\widetilde{S}_n)$. Let $V_{n,m}=\widetilde{S}_n\cap
q^{-1}(S_{n,m})$ for every $n,m\in \N$. It is easy to see that
$V_{n,m}$ are relatively weakly open. We will now prove that they
are the sets we need.

Let $B\subseteq A$ be convex and such that $B\cap
V_{n,m}\neq\emptyset$ for all $n,m\in\N$. Fix some
$\varepsilon>0$, and denote $B_{\varepsilon}=\{x\in A\,: \,
\dist(x,B)<\varepsilon\}$. Evidently, $B_{\varepsilon}$ is an open
convex set intersecting all the $V_{n,m}$. Fixed $n\in\N$, we have
that
$$
B_{\varepsilon}\cap V_{n,m}=B_{\varepsilon}\cap
\widetilde{S}_n\cap q^{-1}(S_{n,m})\neq\emptyset,
$$
so
$$
q\bigl(B_{\varepsilon}\cap\widetilde{S}_n\bigr)\cap S_{n,m}\neq
\emptyset
$$
and the choice of $S_{n,m}$ allows us to get that
$$
\cconv\bigl(q(B_{\varepsilon}\cap\widetilde{S}_n)\bigr)=
\overline{q\bigl(B_{\varepsilon}\cap\widetilde{S}_n\bigr)} \supseteq
q(\widetilde{S}_n).
$$
Notice that $B_{\varepsilon}\cap\widetilde{S}_n$ is open and
convex, hence, so is $q(B_{\varepsilon}\cap\widetilde{S}_n)$. This
implies that the interior of the set
$\overline{q(B_{\varepsilon}\cap\widetilde{S}_n)}$ coincides with
$q(B_{\varepsilon}\cap\widetilde{S}_n)$. Now, using that
$q(\widetilde{S}_n)$ is open, we get that
$$
q\bigl(B_{\varepsilon}\cap\widetilde{S}_n\bigr)\supseteq
q(\widetilde{S}_n)
$$
and, in particular, $q(B_{\varepsilon}\cap\widetilde{S}_n)\ni
q(a)$. This means that there exists $x_n\in
B_{\varepsilon}\cap\widetilde{S}_n$, such that $q(x_n)=q(a)$,
i.e.\ that $x_n\in B_{\varepsilon}\cap S_n$. Since
$B_{\varepsilon}\subseteq A$ and $\{S_n\}$ is a determining
sequence for $A_a$, we get that $B_{\varepsilon}\supseteq A_a$.
Finally, the arbitrariness of $\varepsilon$ implies that
$\overline{B}\supseteq A_a\ni a$.
\end{proof}

Let us state two immediate consequences of this result.

\begin{corollary}
Let $X$ be a separable Banach space which is not SCD.
\begin{enumerate}
\item[(a)] $X$ contains copies of $\ell_1$, and the quotient
    of $X$ over any copy of $\ell_1$ also contains $\ell_1$.
\item[(b)] Consequently, for every $\ell_1$ subspace $Y_1$ of
    $X$, there is another $\ell_1$ subspace $Y_2$ such that
    $Y_1$ and $Y_2$ are mutually complemented in the closed
    linear span of $Y_1+Y_2$ (i.e.\ $\overline{Y_1+Y_2}=
    Y_1+Y_2 = Y_1 \oplus Y_2$). In particular, $Y_1 \cap
    Y_2=0$.
\end{enumerate}
\end{corollary}

\begin{proof}
(a) is immediate from the above theorem and
Theorem~\ref{Examples-Theorem:ell1}. (b) follows from (a) and the
``lifting'' property of $\ell_1$ \cite[Proposition~2.f.7]{L-T}.
\end{proof}

One may wonder whether item (b) of the above corollary can
actually be a characterization of those separable Banach spaces
which are not SCD. This is not the case as the following remark
shows.

\begin{remark}
{\slshape The space $X= \ell_2(\ell_1)$ (which is an SCD space,
even more it has the RNP) has the following property: it contains
isomorphic copies of $\ell_1$ and for every $\ell_1$ subspace $Y
\subseteq X$, there is another $\ell_1$ subspace $Z \subseteq X$,
such that $Z$ and $Y$ are mutually complemented in the closed
linear span of $Y+Z$.\ }
\end{remark}

\begin{proof}
Let $\{X_n\}_{n=1}^{\infty}$ be a sequence of isometric copies of
$\ell_1$. Then, $X$ is isometric to the $\ell_2$ direct sum of the
spaces $X_n$, $\left[\bigoplus_{n\in\N} X_n\right]_{\ell_2}$. Fix
an $\ell_1$-subspace $Y \subseteq X$ and let us prove that some of
the $X_n$ can be taken as $Z$. Assume to the contrary that for
every $n \in \N$
$$
\inf\{\|y - x\| : y \in S_Y, x \in X_n\}=0.
$$
Then, for every $n \in \N$ there are $y_n \in S_Y$ and $x_n \in
X_n$ with $\|y_n - x_n\| < 10^{-n}$. Since $(x_n)$ forms a bounded
sequence of disjoint elements, $(x_n) \longrightarrow 0$ in the
weak topology. But then $(y_n) \longrightarrow 0$ in the weak
topology as well, which is impossible since $(y_n) \subseteq S_Y$
and $Y$ has the Schur property.
\end{proof}

\begin{corollary}\label{cor:finite-sums}
Let $X_1,\ldots,X_n$ be SCD Banach spaces. Then, $X_1\oplus
\cdots\oplus X_n$ is SCD.
\end{corollary}

Our next goal is to deal with infinite sums. To do so, we need to
recall the concept of unconditional sums. Given a sequence
$\left\{(X_n,\|\cdot\|_n)\,:\,n\in \N\right\}$ of Banach spaces,
and a Banach space $E$ of sequences whose norm satisfies
$$
\|(t_i)\|_E=\|(|t_i|)\|_E \qquad \bigl((t_i)\in E\bigr),
$$
we denote by $\left[\bigoplus_{n\in\N} X_n\right]_{E}$ the Banach
space of all sequences $(x_n)\in\prod_{n=1}^{\infty} X_n$, so that
$$
\|(x_n)\|=\|(\|x_n\|_n)\|_E<\infty.
$$

\begin{theorem}\label{teo1}
Let $\{X_n\,:\,n\in\N\}$ be a sequence of SCD spaces and let $E$
be a Banach space of sequences whose canonical basis is a
$1$-unconditional and shrinking basis (i.e.\ $E$ does not contain
copies of $\ell_1$). Then, $X= \left[\bigoplus_{n\in\N}
X_n\right]_{E}$ is also an SCD space.
\end{theorem}

\begin{proof}
For every $m\in \N$, we denote
$$
Y_m=\bigl[X_1\oplus X_2\oplus \ldots\oplus X_m\oplus 0 \oplus 0
\oplus \cdots\bigr]_E \subseteq X
$$
and let $P_m:X\longrightarrow Y_m$ be the natural projection. Let
$A$ be a convex bounded subset of $X$. Now, for every $m\in
\mathbb{N}$, $P_m(A)$ is a convex bounded subset of $Y_m$, which
is an SCD space by Corollary~\ref{cor:finite-sums}. Hence, there
is a determining sequence $\{S_{m,k}\,:\, k\in\N\}$ of slices of
$P_m(A)$. Consider $\widetilde{S}_{m,k} =
P_m^{-1}\bigl({S}_{m,k}\bigr) \cap A$. We will prove that
$\{\widetilde{S}_{m,k}\,:\,k,m\in\N\}$ is a determining countable
collection of slices of $A$.

Let $B$ be a subset of $A$ intersecting all the
$\widetilde{S}_{m,k}$. We fix an arbitrary point $a\in A$ and we
will prove that $a\in \cconv (B)$. Since $B$ intersects all the
$\widetilde{S}_{m,k}$, $P_m(B)$ intersects $S_{m,k}$ for every
integer $k$. It follows that $\cconv \bigl(P_m(B)\bigr)\supseteq
P_m(A)$. In particular, $\cconv \bigl(P_m(B)\bigr)\ni P_m(a)$.
That means that there exists $b_m\in \conv B$ such that
$\|P_m(b_m-a)\|<\frac{1}{m}$. Then, it is easy to see that $b_m$
tends to $a$ coordinate-wise. But since the canonical basis of $E$
is at the same time a shrinking basis, we get that $b_m$ tends to
$a$ in the weak topology. So we can apply Mazur's theorem and get
a sequence $\{b_m'\}$ with $b_m'\in\conv\bigl(\{b_k\,:\,k\geqslant
m \}\bigr) \subseteq\conv (B)$ which tends to $a$ in the norm
topology. But this exactly means that $a\in\cconv (B)$, which was
to be proved.
\end{proof}

The next result deals with unconditional sums when the natural
basis of $E$ is boundedly complete. Its proof, which is more bulky
than the above one, needs a preliminary result which can be of
independent interest.

Let $X$ be a Banach space, $A$ be a convex set in $X$ and
$\varepsilon$ be a positive real. A point $a\in A$ is called an
\emph{$\varepsilon$-accessible point} of $A$ if there is a
sequence $\{V_n\,:\,n\in \N\}$ of relatively weakly-open subsets
of $A$, such that for every $B\subseteq A$, if $B$ intersects all
the $V_n$, then $\dist(a,\conv B)<\varepsilon$.

\begin{lemma}\label{examples-lem2}
Let $X$ be a Banach space and let $A$ be a separable convex
bounded subset of $X$. Suppose that for every convex $C \subseteq
A$ and every $\eps>0$, there is an $\eps$-accessible point in $C$.
Then, $A$ is an SCD set.
\end{lemma}

\begin{proof}
Notice that, since $A$ is separable, to prove this lemma it is
enough to show that for every $\varepsilon>0$, the set
$A_{\varepsilon}$ of $\varepsilon$-accessible points of $A$ is
dense in $A$. Since $A_{\varepsilon}$ is convex, it is enough to
show that $A_{\varepsilon}$ is weakly dense in $A$. Fix some
convex relatively weakly-open subset $V \subseteq A$. By the
assumption, there is an $\varepsilon$-accessible point of $V$. But
this point is also an $\varepsilon$-accessible point of $A$ since
$V$ is relatively weakly-open.
\end{proof}

We are now able to state and prove the second result for
unconditional sums.

\begin{theorem}\label{Examples-th4}
Let $\{X_n\,:\,n\in \N\}$ be a sequence of SCD spaces and let $E$
be a space of sequences whose natural basis is a $1$-unconditional
and boundedly complete basis (i.e.\ $E$ does not contain
isomorphic copies of $c_0$). Then, $X= \left[\bigoplus_{n\in\N}
X_n\right]_{E}$ is an SCD space.
\end{theorem}

\begin{proof}
Let a convex bounded subset $A$ of $X$ and $\varepsilon>0$ be
fixed. Consider the subset
$$
 A_E = \bigl\{(a_n)_{n\in \N} \in E\, : \, \exists x
=(x_n)_{n\in \N} \in A \text{ with } \|x_n\| = |a_n| \text{ for
all } n\in \N\}.
 $$
Since $A_E$ is a bounded subset of a space with the RNP, there are
a functional $b =(b_n)_{n\in \N} \in E^*$ and a positive number
$\alpha$ such that the slice
$$
S(A_E) = \bigl\{(a_n)_{n\in
\N} \in A_E\ : \ \sum_{n\in \N} b_n a_n > \alpha\bigr\}
$$
has diameter smaller than $\eps/4$. Taking into account that $A_E$
is symmetric, we may assume that $b_n \geq 0$ (the slice of $A_E$
defined by $|b| =(|b_n|)_{n\in \N}$ is isometric to $S(A_E)$). Fix
an $x \in A$ with $(\|x_n\|)_{n\in \N}\in S(A_E)$ and pick $x_n^*
\in S_{X_n^*}$ such that $x_n^*(x_n)=\|x_n\|$. Write $f_n= b_n
x_n^*$, $f = (f_n)_{n\in \N} \in X^*$. We claim that for the slice
$$
S = \left\{(x_n)_{n\in \N} \in
A \,:\, \sum_{n\in \N} f_n(x_n) > \alpha\right\}
$$
there is an $m\in \N$ with the following property
\begin{equation}\label{eq:PROPERTY}
\bigl\|(0, \ldots, 0, y_{m+1}, y_{m+2}, \ldots)\bigr\|
<\frac{\eps}{2} \qquad \text{for all } (y_n)_{n\in \N} \in S.
\end{equation}
To show this, it is sufficient to select $m$ in such a way that
$\|(0, \ldots, 0, x_{m+1}, x_{m+2}, \ldots)\| < \eps/4$ and to use
that $\diam~S(A_E) < \eps/4$. In fact, with such a choice of $m$
we get
\begin{align*}
\|(0, \ldots, 0, y_{m+1}, y_{m+2}, \ldots)\| &= \|(0, \ldots, 0,
\|y_{m+1}\|, \|y_{m+2}\|, \ldots)\| \\
&\leq \bigl\|(0, \ldots, 0, \|x_{m+1}\|, \|x_{m+2}\|, \ldots)\bigr\|
 + \\ &
\qquad \quad + \bigl\|(0, \ldots, 0, \bigl| \|x_{m+1}\| -
\|y_{m+1}\| \bigr|, \bigl| \|x_{m+2}\| - \|y_{m+2}\| \bigr|,
\ldots)\bigr\| \\ &\leq \frac\eps4 + \left\|(\bigl| \|x_{1}\| - \|y_{1}\| \bigr|, \bigl| \|x_{2}\| -
\|y_{2}\| \bigr|, \ldots)\right\| \leq \frac\eps2.
\end{align*}
Let us prove that $x$ is an $\varepsilon$-accessible point of $A$.
Consider
$$
Y_m=\bigl[X_1\oplus X_2\oplus \ldots\oplus X_m\oplus 0 \oplus 0
\oplus \cdots\bigr]_E \subseteq X
$$
and $P_m:X\longrightarrow Y_m$ the natural projection. By
Corollary~\ref{cor:finite-sums}, $Y_m$ is an SCD space and, since
$P_m(S)$ is a convex bounded set in $Y_m$, there exists a
determining sequence $\{S_n\,:\,n\in\N\}$ of slices of $P_m(S)$.
Notice that $Y_m^*$ isometrically embeds into $X^*$. For every
integer $n\in\N$, we consider $\widetilde{S}_n= P_m^{-1} {S}_n
\cap S$, which is a slice of $S$ and, obviously, relatively
weakly-open in $A$. Let $B$ be a subset of $A$ which intersects
all the $\widetilde{S}_n$. We'll now prove that then
$\dist\bigl(x,\conv(B)\bigr)<\varepsilon$. Since $B$ intersects
all the $\widetilde{S}_n$, we can find a sequence
$\{y_n\}\subseteq B$, such that $y_n\in \widetilde{S}_n$ for every
$n\in\N$. This implies that $P_m(y_n)\in S_n$ for all $n \in \N$
and so $\cconv\bigl(\{P_m(y_n)\,:\,n\in\N\}\bigr) \supseteq
P_m(S)$. In particular, $P_m(x)\in
\cconv\bigl(\{P_m(y_n)\,:\,n\in\N\}\bigr)$. But
\eqref{eq:PROPERTY} gives us that the $m$-th tails of $x$ and of
all the $y_n$ are small, that is,
$$
\|x-P_m(x)\|<\frac{\varepsilon}{2} \quad \text{and} \quad \|y_n - P_m(y_n)\|
< \eps/2 \qquad \bigl(\text{for all}\ n\in \N\bigr).
$$
This gives us that $\dist\bigl(a,\conv(B)\bigr)<\varepsilon$ and
the proof is complete.
\end{proof}

An immediate consequence is the following.

\begin{example}
{\slshape The spaces $c_0(\ell_1)$ and $\ell_1(c_0)$ are SCD.\ }
\end{example}

This result, together with those results of
section~\ref{sec:sets}, gives us the following examples.

\begin{example}
{\slshape The spaces $c_0 \otimes_\eps c_0$, $c_0 \otimes_\pi
c_0$, $c_0 \otimes_\eps \ell_1$, $c_0 \otimes_\pi \ell_1$,
$\ell_1\otimes_\eps \ell_1$, and $\ell_1 \otimes_\pi \ell_1$ are
SCD.\ } Indeed, it is well known that $c_0\otimes_\eps c_0\equiv
c_0$, $c_0\otimes_\eps \ell_1\equiv c_0(\ell_1)$, $c_0 \otimes_\pi
\ell_1\equiv \ell_1(c_0)$, and $\ell_1 \otimes_\pi \ell_1\equiv
\ell_1$ (see \cite[Examples 2.19 and 3.3]{Ryan}, for instance), so
these cases are clear from the above example. For the remaining
cases, just observe that $\bigl[c_0\otimes_\pi c_0\bigr]^*\equiv
\ell_1\otimes_\eps \ell_1$ (since $[c_0\otimes_\pi c_0]^*\equiv
L(c_0,\ell_1)$ \cite[p.~24]{Ryan}, $K(c_0,\ell_1)\equiv
\ell_1\otimes_\eps \ell_1$ \cite[Corollary~4.13]{Ryan} and
$K(c_0,\ell_1)=L(c_0,\ell_1)$ since $\ell_1$ has the Schur
property and $c_0^*$ is separable), so $c_0\otimes_\pi c_0$ is
Asplund and $\ell_1\otimes_\eps \ell_1$ has the RNP.
\end{example}

Since for $X$ and $Y$ being $c_0$ or $\ell_1$ one has
$K(X,Y)\equiv X^*\otimes_\eps Y$ \cite[Corollary~4.13]{Ryan}, the
following examples follow.

\begin{example}
{\slshape The spaces $K(c_0)$ and $K(c_0,\ell_1)$ are SCD. The
spaces $K(\ell_1)$ and $K(\ell_1,c_0)$ contain $\ell_\infty$ and
so they are not separable, all the more not SCD.\ }
\end{example}

Another example in this line is the following.

\begin{example}
{\slshape The spaces $\ell_2\otimes_\pi \ell_2 \equiv
\mathcal{L}_1(\ell_2)$, and $\ell_2\otimes_\eps \ell_2\equiv
K(\ell_2)$ are SCD.\ } Indeed, the first space has the RNP and the
second is an Asplund space.
\end{example}

\section{An application to spaces with numerical index~$1$}
\label{sec:applications-to-ni1} Our aim in this section is to show
that SCD spaces with the alternative Daugavet property are lush.
To get such a result, we need to establish a characterization of
the alternative Daugavet property which can be of independent
interest. We first recall a previous characterization in terms of
slices.

\begin{lemma}[\textrm{\cite[Proposition~2.1]{MaOi}}] \label{lemma:ADP}
A Banach space $X$ has the alternative Daugavet property if and
only if for every $x \in S_X$, every $\eps > 0$ and every slice
$S$ of $B_X$, there is a $y \in S$ such that $\max_{\theta \in
\T}\|x + \theta y\| > 2 - \eps$.
\end{lemma}

We need some notation. Denote $K(X^*)$ the weak$^*$-closure in
$X^*$ of $\ext[X^*]$, and for every slice $S$ of $B_X$ and every
$\eps > 0$, we write
\begin{align*}
D(S,\eps)&=\bigl\{y^* \in K(X^*)\, :\, S \cap \T S(B_X,y^*,\eps) \neq
\emptyset\bigr\} \\ &=\bigl\{y^* \in K(X^*)\, :\, S \cap \caconv\bigl(S(B_X,y^*,\eps)\bigr) \neq
\emptyset \bigr\},
\end{align*}
which is relatively weak$^*$-open in $K(X^*)$. Here is the
promised characterization of the alternative Daugavet property.

\begin{prop} \label{adp-prop1}
For a Banach space $X$, the following assertions are equivalent:
\begin{enumerate}
\item[(i)] $X$ has the alternative Daugavet property.
\item[(ii)] For every $x \in S_X$, every $\eps > 0$ and every
    slice $S \subseteq B_X$,  there is $y^* \in K(X^*)$ such
    that $x \in S(B_X,y^*,\eps)$ and $S \cap \T
    S(B_X,y^*,\eps) \neq \emptyset$.
\item[(iii)] For every $x \in S_X$, every $\eps > 0$ and every
    slice $S \subseteq B_X$, there is $y^* \in D(S,\eps)$ such
    that $x \in S(B_X,y^*,\eps)$.
\item[(iv)] For every $\eps > 0$ and every slice $S \subseteq
    B_X$, the set $D(S,\eps)$ is weak$^*$-dense in $K(X^*)$.
\item[(v)] For every $\eps > 0$ and every sequence
    $\{S_n\,:\,n\in\N\}$ of slices of $B_X$, the set
    $\bigcap_{n \in \N}D(S_n,\eps)$ is weak$^*$-dense in
    $K(X^*)$.
\end{enumerate}
\end{prop}

\begin{proof} The implications (i) $\Longleftrightarrow$ (ii)
$\Longleftrightarrow$ (iii) are easy consequences of
Lemma~\ref{lemma:ADP}.

(iii) $\Longrightarrow$ (iv). To show weak$^*$-density of
$D(S,\eps)$  in $K(X^*)$ it is sufficient to demonstrate that the
weak$^*$ closure of $D(S,\eps)$ contains every extreme point $x^*$
of $S_{X^*}$. Since weak$^*$-slices form a base of neighborhoods
of $x^*$ in $B_{X^*}$, it is sufficient to prove that every
weak$^*$-slice $S(B_{X^*},x,\delta)$ with $\delta \in (0, \eps)$
intersects $D(S,\eps)$, i.e.\ that there is a point $y^* \in
D(S,\eps)$, such that $y^* \in S(B_{X^*},x,\delta)$. But we know
that there is a point $y^* \in D(S,\delta) \subseteq D(S,\eps)$,
such that $x \in S(B_{X},y^*, \delta)$, which means that $y^* \in
S(B_{X^*},x,\delta)$.

(iv) $\Longrightarrow$ (iii). If $D(S,\eps)$ is weak$^*$-dense in
$K(X^*)$, then for every $x \in S_X$ there is a $y^* \in
D(S,\eps)$ such that $x \in S(B_X,y^*,\eps)$.

The remaining equivalence (iv) $\Longleftrightarrow$ (v) follows
from the fact that $D(S,\eps)$ is not only weak$^*$-dense but also
weak$^*$-open, and $K(X^*)$ is weak$^*$-compact, so Baire's
theorem is applicable.
\end{proof}

It is possible to give a result analogous to the above one for the
Daugavet property. We need to change a little bit the notation.
For every slice $S$ of $B_X$ and every $\eps > 0$, we write
\begin{align*}
\widetilde{D}(S,\eps)&=\bigl\{y^* \in K(X^*)\, :\, S \cap S(B_X,y^*,\eps) \neq
\emptyset\bigr\} \\ &=\bigl\{y^* \in K(X^*)\, :\, S \cap \cconv\bigl(S(B_X,y^*,\eps)\bigr) \neq
\emptyset \bigr\}
\end{align*}
which is relatively weak$^*$-open in $K(X^*)$. The proof of the
next result is almost the same as the above one, replacing
Lemma~\ref{lemma:ADP} by \cite[Lemma~2.2]{KadSSW}. We include it
here for future use.

\begin{prop}\label{prop:charac-DPr}
For a Banach space $X$, the following assertions are equivalent:
\begin{enumerate}
\item[(i)] $X$ has the Daugavet property.
\item[(ii)] For every $x \in S_X$, every $\eps > 0$ and every
    slice $S \subseteq B_X$,  there is $y^* \in K(X^*)$ such
    that $x \in S(B_X,y^*,\eps)$ and $S \cap S(B_X,y^*,\eps)
    \neq \emptyset$.
\item[(iii)] For every $x \in S_X$, every $\eps > 0$ and every
    slice $S \subseteq B_X$, there is $y^* \in
    \widetilde{D}(S,\eps)$ such that $x \in S(B_X,y^*,\eps)$.
\item[(iv)] For every $\eps > 0$ and every slice $S \subseteq
    B_X$, the set $\widetilde{D}(S,\eps)$ is weak$^*$-dense in
    $K(X^*)$.
\item[(v)] For every $\eps > 0$ and every sequence
    $\{S_n\,:\,n\in\N\}$ of slices of $B_X$, the set
    $\bigcap_{n \in \N}\widetilde{D}(S_n,\eps)$ is
    weak$^*$-dense in $K(X^*)$.
\end{enumerate}
\end{prop}

We are now ready to show the main result of this section.

\begin{theorem} \label{ssd-th1}
Every Banach space $X$ with the alternative Daugavet property
whose unit ball is an SCD set is lush. In particular, every SCD
space with the alternative Daugavet property is lush.
\end{theorem}

\begin{proof}
Let $\{S_n\,:\,n\in\N\}$ be the sequence of slices of $B_X$ from
the definition of an SCD set. Then, by
Proposition~\ref{adp-prop1}.v, for every $\eps > 0$ the set
$\bigcap_{n \in \N} D(S_n,\eps)$ is weak$^*$-dense in $K(X^*)$.
So, for every $x \in S_X$ there is $y^* \in \bigcap_{n \in
\N}D(S_n,\eps)$ such that $x \in S(B_X,y^*,\eps)$. According to
the definition of $D(S_n,\eps)$, this means that $S_n \cap \caconv
\bigl(S(B_X,y^*,\eps)\bigr) \neq \emptyset$ for all $n \in \N$.
Then, we obtain that $\caconv \bigl(S(B_X,y^*,\eps)\bigr) = B_X$,
which implies lushness of $X$ \cite[Theorem~2.1]{BKMM-lush}.
\end{proof}

\begin{remark}\label{remark:ADP=lush_with_almostSCD}
Let us observe that in the above proof a (formally) weaker version
of an SCD set is used. A convex bounded subset $A$ of a Banach
space $X$ is said to be \emph{almost slicely countably determined}
(\emph{almost-SCD} in short) if there is a sequence
$\{V_n\,:\,n\in\N\}$ of subsets of $A$ such that for every
$B\subseteq A$ intersecting all the $V_n$, one has $\caconv (B)
\supseteq A$. The proof of the above theorem actually shows that
{\slshape every Banach space $X$ with the alternative Daugavet
property whose unit ball is an almost-SCD is lush.\ }
\end{remark}

Theorem~\ref{ssd-th1} has already been known for Asplund spaces
and for spaces with the RNP \cite[Remark~6]{LMP1999}, regardless
of the separability (necessary for the SCD and so for our result).
Our next goal is to particularize Theorem~\ref{ssd-th1} to more
cases where we are able to remove the separability. The proof of
the following results is a consequence of the facts that lushness
and the alternative Daugavet property are separably determined
(see \cite[Theorem~4.2]{BKMM-lush} for the first case and the
remark below for the second one).

\begin{remark}\label{rem:ADP-separably-determined}
It is shown in \cite[Theorem~4.5]{KadSW2} that the Daugavet
property is separably determined. With a little effort, the proof
can be adapted to the alternative Daugavet property: {\slshape A
Banach space $X$ has the alternative Daugavet property if and only
if for every separable subspace $Y \subseteq X$ there is a
separable subspace $Z \subseteq X$ which contains $Y$ and has the
alternative Daugavet property.\ }
\end{remark}

\begin{corollary}\label{cor:stronglyregular-ADP=>lush}
Let $X$ be a Banach space with the alternative Daugavet property.
If $X$ is strongly regular (in particular, CPCP), then $X$ is
lush.
\end{corollary}

\begin{corollary}\label{cor:nonell1-ADP=>lush}
Let $X$ be a Banach space with the alternative Daugavet property.
If $X$ does not contain $\ell_1$, then $X$ is lush.
\end{corollary}

This latter result solves in the positive Problem~$32$ of
\cite{KaMaPa} and it can be used to prove a necessary isomorphic
condition for a real Banach space to have the alternative Daugavet
property.

\begin{corollary}
let $X$ be an infinite-dimensional real Banach space with the
alternative Daugavet property. Then, $X^*$ contains $\ell_1$.
\end{corollary}

\begin{proof}
If $X$ contains $\ell_1$, then $X^*$ contains a quotient
isomorphic to $\ell_\infty$, so $X^*$ contains $\ell_1$ as a
quotient and the ``lifting'' property of $\ell_1$
\cite[Proposition~2.f.7]{L-T} gives us $X^*\supseteq \ell_1$.
Otherwise, Corollary~\ref{cor:nonell1-ADP=>lush} gives us that $X$
is lush. But the dual of an infinite-dimensional real lush space
contains $\ell_1$ \cite[Corollary~4.8]{KMMP}.
\end{proof}

In particular, since Banach spaces with numerical index~$1$ have
the alternative Daugavet property, we get the following corollary
which answers in the positive Problem~$18$ of \cite{KaMaPa}.

\begin{corollary}\label{cor:answer-problem18}
Let $X$ be an infinite-dimensional real Banach space with
$n(X)=1$. Then, $X^*\supseteq \ell_1$.
\end{corollary}

Let us comment that very recently it has been shown that there are
Banach spaces with numerical index~$1$ which are not lush
\cite{KMMS}, so the above result is not covered by
\cite[Corollary~4.9]{KMMP}.

\section{SCD operators}\label{sec:SCD-operators}
\begin{definition}\label{op-def1}
Let $X$ and $Y$ be Banach spaces. A bounded linear operator
$T:~X\longrightarrow Y$ is said to be an \emph{SCD-operator} if
$T(B_X)$ is an SCD set.
\end{definition}

By just recalling the examples of SCD sets and SCD spaces given in
sections \ref{sec:sets} and \ref{sec:spaces}, we get the main
examples of SCD-operators.

\begin{examples} Let $X$ and $Y$ be Banach spaces and let $T:~X\longrightarrow
Y$ be a bounded linear operator such that $T(X)$ is separable.
\begin{enumerate}
\item[(a)] {\slshape If $T(B_X)$ has small combinations of
    slices, then $T$ is an SCD-operator.\ }
\item[(b)] In particular, {\slshape if $T(B_X)$ is a
    Radon-Nikod\'{y}m set (i.e.\ if $T$ is a strong Radon-Nikod\'{y}m
    operator), then $T$ is an SCD-operator.\ }
\item[(c)] {\slshape If $T(B_X)$ does not contain
    $\ell_1$-sequences, then $T$ is an SCD-operator.\ }
\item[(d)] In particular, {\slshape if $T$ does not fix copies
    of $\ell_1$, then $T$ is an SCD-operator.\ } Indeed, if
    $T(B_X)$ contains an $\ell_1$-sequence $(T e_n)_{n\in\N}$
    with $e_n\in B_X$ ($n\in \N$), then as in the proof of the
    ``lifting'' property of $\ell_1$
    \cite[Proposition~2.f.7]{L-T},
    $Y=\overline{\textrm{lin}}\{e_n\,:\,n\in\N\}$ is a copy of
    $\ell_1$ and $T|_Y$ is an isomorphic embedding, a
    contradiction (see \cite[Proposition~1]{Weis}).
\end{enumerate}
\end{examples}

The aim of this section is to show that SCD-operators behave in a
very good way with respect to the Daugavet and the alternative
Daugavet equations. We start with the best result we can get for
the alternative Daugavet property.

\begin{theorem}\label{theorem:ADP-SCD-opertor}
Let $X$ be a Banach space with the alternative Daugavet property
and let $T\in L(X)$ be an SCD-operator. Then,
$\max\limits_{\theta\in\T}\|\Id+\theta\, T\|=1+\|T\|$.
\end{theorem}

\begin{proof}
Without loss of generality, we may assume that $\|T\|=1$. We take
a determining sequence $\{S_n\,:\,n\in\N\}$ of slices of $T(B_X)$
and we notice that the sets $T^{-1}(S_n)\cap B_X$ are slices of
$B_X$. Given $\eps>0$ fixed, we take $a\in S_X$ such that
$\|T(a)\|>1-\eps$. Now, Proposition~\ref{adp-prop1}.v gives us
that $\bigcap_{n\in\N}D\bigl(T^{-1}(S_n),\varepsilon\bigr)$ is
weak$^{*}$-dense in $K(X^*)$ (which is norming for $X$), so we may
find $y^*\in\bigcap_{n\in\mathbb{N}}
D\bigl(T^{-1}(S_n),\varepsilon\bigr)$ such that
\begin{equation}\label{eq:ADP-SCD-1}
\re y^{*}(T(a))\geq \|T(a)\|-\varepsilon>1 - 2\eps.
\end{equation}
By the definition of $D\bigl(T^{-1}(S_n),\varepsilon\bigr)$, we
get that
$$
\caconv \bigl(S(B_X,y^*,\varepsilon)\bigr)\cap
T^{-1}(S_n)\neq\emptyset \qquad (n\in \N).
$$
Thus, $T\bigl(\caconv\bigl(S(B_X,y^*,\varepsilon)\bigr)\bigr) \cap
S_n\neq\emptyset$ for all $n\in \N$, and using the fact that
$\{S_n\,:\,n\in\N\}$ is determining, we deduce that
$$
\overline{T\bigl(\aconv \bigl(S(B_X,y^*,\varepsilon)\bigr)\bigr)}=
\caconv\bigl(T\bigl(\caconv\bigl(S(B_X,y^*,\varepsilon)\bigr)\bigr)\supseteq
T(B_X).
$$
In particular, $T(a)\in \overline{T\bigl(\aconv
\bigl(S(B_X,y^*,\varepsilon)\bigr)\bigr)}$, which means that there
is
$$
z\in T\bigl(\aconv\bigl(S(B_X,y^*,\varepsilon)\bigr)\bigr)
\quad \text{with} \quad \|T(a)-z\|<\varepsilon,
$$
and it follows from \eqref{eq:ADP-SCD-1} that
\begin{equation}\label{eq:ADP-SCD-2}
\re y^*(z)>1-3\varepsilon.
\end{equation}
Notice that $z$ can be represented in the following way
$$
z=T\left(\sum\limits_{k=1}^m\lambda_k \theta_k
x_k\right)=\sum\limits_{k=1}^m\lambda_k \,\theta_k T(x_k)
$$
where $x_k\in S(B_X,y^*,\varepsilon)$, $\theta_k\in\T$,
$\lambda_k\geq 0$ for $k=1,\ldots,m$ and
$\sum_{k=1}^m\lambda_k=1$. Then, it follows from
\eqref{eq:ADP-SCD-2} that there exists $k_0\in \{1,\ldots,m\}$
such that
$$
\re y^*\bigl(\theta_{k_0}\,T(x_{k_0})\bigr)>1-3\varepsilon.
$$
Now, since $x_{k_0}\in S(B_X,y^*,\varepsilon)$, we get that
$$
\re y^*\bigl(x_{k_0}+\theta_{k_0}\,T(x_{k_0})\bigr)>2-4\varepsilon.
$$
It follows that
$$
\|\Id+\theta_{k_0}T\|\geq \|x_{k_0}+\theta_{k_0}\,T(x_{k_0})\| \geq
\re y^*\bigl(x_{k_0}+\theta_{k_0}\,T(x_{k_0})\bigr)>2-4\varepsilon.
$$
Finally, the arbitrariness of $\varepsilon$ gives the result.
\end{proof}

\begin{remark}\label{remark:ADPalmostSCD-operators}
Analogously to the situation described in
Remark~\ref{remark:ADP=lush_with_almostSCD}, in the above proof we
have used a formally weaker property than being an SCD-operator.
Therefore, the result proved is the following. {\slshape Let $X$
be a Banach space with the alternative Daugavet property and let
$T\in L(X)$ such that $T(B_X)$ is an almost-SCD set. Then,
$\max\limits_{\theta\in\T}\|\Id+\theta\, T\|=1+\|T\|$.\ }
\end{remark}

We can easily obtain a version of
Theorem~\ref{theorem:ADP-SCD-opertor} for operators with non
separable range which is useful for applications.

\begin{corollary}\label{cor:ADP-SCD-operators-nonseparable}
Let $X$ be a Banach space with the alternative Daugavet property
and let $T\in L(X)$ be such that $T(B_Y)$ is an SCD set for every
separable subspace $Y$ of $X$. Then,
$\max\limits_{\theta\in\T}\|\Id+\theta\, T\|=1+\|T\|$.
\end{corollary}

\begin{proof}
We first take a separable subspace $Y_1$ of $X$ such that
$\|T|_{Y_1}\|=\|T\|$. Then,
Remark~\ref{rem:ADP-separably-determined} provides us with a
separable subspace $Y_2$ with the alternative Daugavet property
which contains $\bigcup_{k=0}^\infty T^k(Y_1)$. We apply again
Remark~\ref{rem:ADP-separably-determined} to get a separable
subspace $Y_3$ with the alternative Daugavet property which
contains $\bigcup_{k=0}^\infty T^k(Y_2)$, and so on. Then, the
space $Y=\overline{\bigcup_{n\in\N} Y_n}$ is separable,
$T$-invariant, $\|T|_Y\|=\|T\|$, and it has the alternative
Daugavet property (just use Lemma~\ref{lemma:ADP}). Since $T(B_Y)$
is SCD, Theorem~\ref{theorem:ADP-SCD-opertor} gives us that
\begin{equation*}
\max\limits_{\theta\in\T}\|\Id + \theta\,T\|\geqslant \max\limits_{\theta\in\T}
\left\|\Id|_Y+\theta\, T|_Y\right\|=1+\left\|T|_Y\right\|=1+\|T\|.\qedhere
\end{equation*}
\end{proof}

The following particular cases are especially interesting. The
first one solves \cite[Problem~$33$]{KaMaPa}.

\begin{corollary}
Let $X$ be a Banach space with the alternative Daugavet property
and let $T\in L(X)$ be an operator which does not fix copies of
$\ell_1$. Then, $\max\limits_{\theta\in\T}\|\Id+\theta\,
T\|=1+\|T\|$.
\end{corollary}

\begin{corollary}
Let $X$ be a Banach space with the alternative Daugavet property
and let $T\in L(X)$ be an operator such that $T(B_X)$ is strongly
regular. Then, $\max\limits_{\theta\in\T}\|\Id+\theta\,
T\|=1+\|T\|$.
\end{corollary}

It is possible to show an analogous result to
Theorem~\ref{theorem:ADP-SCD-opertor} for spaces with the Daugavet
property and the Daugavet equation. Actually, it is possible to
get a better result. We need some notation and preliminary
results. A bounded linear operator $T:X\longrightarrow Y$ between
two Banach spaces $X$ and $Y$ is said to be a \emph{strong
Daugavet operator} if for every $x,y\in S_X$ and every $\eps>0$,
there is an element $z\in S_X$ such that
$$
\|x+z\|\geq 2-\eps \qquad \text{and} \qquad \|Ty - Tz\|<\eps
$$
(see \cite[\S3]{KadSW2} for the definition and the following
properties). If $T\in L(X)$ is a strong Daugavet operator and $X$
has the Daugavet property, then $T$ satisfies Daugavet equation.
On the other hand, finite-rank operators from a space with the
Daugavet property are strong Daugavet operators. Our next goal is
to show that actually, SCD-operators are strong Daugavet
operators.

\begin{prop}\label{prop:appl2-lem2}
Let $X$ be a Banach space with Daugavet property, $Y$ a Banach
space, and let $T:X\longrightarrow Y$ be an SCD-operator. Then,
$T$ is a strong Daugavet operator.
\end{prop}

\begin{proof}
Since $T$ is an SCD-operator, we may find a determining sequence
$\{S_n\,:\,n\in\N\}$ of slices of $T(B_X)$, and we notice that the
sets $T^{-1}(S_n)\cap B_X$ are slices of $B_X$. We fix
$\varepsilon>0$ and $x,y \in S_X$.

Since $X$ has the Daugavet property,
Proposition~\ref{prop:charac-DPr}.v gives us that
$\bigcap_{n\in\N}\widetilde{D}\bigl(T^{-1}(S_n),\frac{\varepsilon}{2}\bigr)$
is weak$^{*}$-dense in $K(X^*)$ (which is norming for $X$), so we
may find $y^*\in\bigcap_{n\in\mathbb{N}}
\widetilde{D}\bigl(T^{-1}(S_n),\frac{\varepsilon}{2}\bigr)$ such
that
\begin{equation}\label{eq:DPr-SCD-1}
x\in S(B_X,y^*,{\textstyle \frac{\varepsilon}{2} }).
\end{equation}
Then, by the definition of
$\widetilde{D}(T^{-1}(S_n),\frac{\varepsilon}{2})$, we have that
$\overline{S(B_X,y^*,{\textstyle \frac{\varepsilon}{2}})}\cap
T^{-1}(S_n)\neq\emptyset$ for every $n\in \N$. Thus,
$$
T\left(\overline{S(B_X,y^*,{\textstyle \frac{\varepsilon}{2}})}\right)\cap
S_n\neq\emptyset \qquad (n\in\N).
$$
Now, since the sequence $\{S_n\,:\,n\in\N\}$ is determining, we
deduce that
$$
T(B_X)\subseteq \cconv\, T\left(\overline{S(B_X,y^*, {\textstyle
\frac{\varepsilon}{2} })}\right)
=\overline{T\left(S(B_X,y^*,{\textstyle \frac{\varepsilon}{2}
})\right)}.
$$
In particular, $Ty\in
\overline{T\left(S(B_X,y^*,\frac{\varepsilon}{2})\right)}$, which
means that there is a $z\in S(B_X,y^*,\frac{\varepsilon}{2})$ such
that
$$
\|Ty-Tz\|<\varepsilon.
$$
Since $x\in S(B_X,y^*,\frac{\varepsilon}{2})$ by
\eqref{eq:DPr-SCD-1}, we also have that
\begin{equation*}
\|x+z\|>2-\varepsilon.
\end{equation*}
Hence, this $z$ meets all the requirements.
\end{proof}

In particular, we obtain the following analogue to
Theorem~\ref{theorem:ADP-SCD-opertor}.

\begin{corollary}
Let $X$ be a Banach space with the Daugavet property. If $T\in
L(X)$ is an SCD-operator, then $\|\Id + T\|=1 + \|T\|$.
\end{corollary}

Our final goal in this section is to get a better result than
Proposition~\ref{prop:appl2-lem2} for a class of operators more
restrictive than the SCD-operators. We need some notation. A
bounded linear operator $T:X\longrightarrow Y$ between two Banach
spaces $X$ and $Y$ is said to be a \emph{narrow operator} if for
every $x^*\in X^*$, the operator
$$
T\tilde{+}\re x^*:X\longrightarrow Y\oplus_1 \R,\qquad  x\longmapsto \bigl(Tx,\re x^*(x)\bigr)
$$
is a strong Daugavet operator (see \cite[\S3 and \S4]{KadSW2} for
this definition and the following properties). Equivalently, $T$
is narrow if and only if for every $x,y\in S_X$, every $\eps>0$,
and every slice $S$ of $B_X$ containing $y$, there is an element
$z\in S$ such that
\begin{equation*}
\|x+z\|\geq 2-\eps \qquad \text{and} \qquad \|Ty - Tz\|<\eps.
\end{equation*}
A narrow operator is strong Daugavet, but the converse result is
not true. It is known that strong Radon-Nikod\'{y}m operators and
operators which do not fix copies of $\ell_1$ from a Banach space
with the Daugavet property are narrow. We are going to extend
these results to the so-called hereditary-SCD-operators.

\begin{definition}\label{appl2-def2}
Let $X$ and $Y$ be Banach spaces. A bounded linear operator
$T:~X\longrightarrow Y$ is said to be a
\emph{hereditary-SCD-operator} if every convex subset of $T(B_X)$
is an SCD set.
\end{definition}

Here is the promised result.

\begin{theorem}\label{appl2-the1}
Let $X$ be a Banach space with Daugavet property and
$T:X\longrightarrow Y$ be a hereditary-SCD-operator. Then, $T$ is
narrow.
\end{theorem}

We need the following lemma, which could be of independent
interest.

\begin{lemma}\label{appl2-lem1}
Let $T:X\longrightarrow Y$ be a hereditary-SCD-operator. Then, for
every $x^*\in X^*$ the operator $T\tilde{+}\re x^*: X
\longrightarrow Y \oplus_1 \R$ is an SCD-operator.
\end{lemma}

\begin{proof}
Denote $P_1:[T\tilde{+}\re x^*](X)\longrightarrow T(X)$ and
$P_2:[T\tilde{+}\re x^*](X)\longrightarrow \R$ the natural
coordinate projections. What we need to show is that there is a
determining sequence of relatively weakly-open subsets of the set
$A=~[T\tilde{+}\re x^*](B_X)$. Since $A$ is separable, it is
enough to prove that for every $a\in A$ there exists a sequence of
relatively weakly open sets $\{V_n\,:\,n\in\N\}$ such that for
every $B\subseteq A$ intersecting all the $V_n$, $a\in\cconv(B)$
(see Remark~\ref{rem:dense-sequence}).

We fix $a\in A$ and denote
$$
A_a=\bigl\{b\in A\,:\, P_1(b)=P_1(a)\bigr\}.
$$
It is easy to see that $A_a$ is of the form
$$
A_a=\bigl\{(P_1(a),t)\,:\, t\in \Delta_a\bigr\},
$$
where $\Delta_a$ is a bounded interval in $\R$. We denote
$\alpha_a=~\inf\Delta_a$ and $\beta_a=\sup\Delta_a$ and we
consider
$$
S_{n,1}=\left\{b\in A\,:\, P_2(b)<\alpha_a+{\textstyle \frac{1}{n}}\right\} \quad \text{and}
\quad S_{n,2}=\left\{b\in A\,:\, P_2(b)>\beta_a-{\textstyle \frac{1}{n}}\right\},
$$
which are non-empty slices of $A$ (since they intersect $A_a$) for
all $n\in\N$ and $i=1,2$. Now, since $T$ is a
hereditary-SCD-operator and $P_1(S_{n,i})\subseteq T(B_X)$ is
convex, for every $n\in\N$ and $i=1,2$ we may find a determining
sequence $\{S_{n,i}^m\,:\,m\in\N\}$ of slices of $P_1(S_{n,i})$.
We write
$$
V_{n,i}^m=S_{n,i}\cap P_1^{-1}(S_{n,i}^m) \qquad (n,m\in\N,\ i=1,2)
$$
which are relatively weakly open subsets of $A$. We will prove
that they are the sets we need. Indeed, let $B\subseteq A$ be such
that
$$
\emptyset\neq B\cap V_{n,i}^m=B\cap S_{n,i}\cap P_1^{-1}(S_{n,i}^m) \qquad (n,m\in\N,\ i=1,2).
$$
For every $n\in\N$, we observe that
$$
P_1(B\cap S_{n,i})\cap S_{n,i}^m\neq \emptyset \qquad (m\in\N,\ i=1,2),
$$
so, since the sequences $\{S_{n,i}^m\,:\,m\in\N\}$ are
determining, we get that
$$
\overline{P_1(B\cap S_{n,i})}=\cconv\bigl(P_1(B\cap S_{n,i})\bigr)
 \supseteq P_1(S_{n,i}) \qquad (i=1,2).
$$
In particular, $\overline{P_1(B\cap S_{n,i})}\ni  P_1(a)$, meaning
that for every $n\in\N$, every $i=1,2$, and every $\eps>0$, there
exists $x_{n,i}^{\eps}\in B\cap S_{n,i}$ such that
\begin{equation}\label{eq:th-narrow-1}
\left\|P_1\left(x_{n,i}^{\varepsilon}\right)-P_1(a)\right\|\leq\eps.
\end{equation}
Now, we fix some $\varepsilon>0$ and, since obviously $a\in A_a$,
we may take $n\in\N$ such that
$$
\alpha_a+\frac{1}{n}-\varepsilon < P_2(a) < \beta_a
-\frac{1}{n}+\varepsilon.
$$
So, for the corresponding $x_{n,1}^{\varepsilon}$ and
$x_{n,2}^{\varepsilon}$, we have
$$
P_2\bigl(x_{n,1}^{\varepsilon}\bigr)-\varepsilon < P_2(a) <
P_2\bigl(x_{n,2}^{\varepsilon}\bigr)+\varepsilon.
$$
Then, there is a convex combination
$$
x_n^{\varepsilon}=\lambda_1\, x_{n,1}^{\varepsilon} + \lambda_2\,
x_{n,2}^{\varepsilon} \qquad (\lambda_1+\lambda_2=1)
$$
(so $x_n^\eps\in \conv(B)$) such that
$$
\bigl|P_2\bigl(x_n^{\varepsilon}\bigr)-P_2(a)\bigr|<\varepsilon.
$$
This, together with \eqref{eq:th-narrow-1}, implies that
$\bigl\|x_n^{\varepsilon}-a\bigr\|<2\varepsilon$, and the
arbitrariness of $\varepsilon>0$ gives us that $a\in \cconv(B)$.
\end{proof}

\begin{proof}[Proof of Theorem~\ref{appl2-the1}]
To prove that $T$ is narrow, it is enough to show that for every
$x^*\in X^*$, the operator $T\tilde{+}\re x^*$ is a strong
Daugavet operator. But this fact follows from
Lemma~\ref{appl2-lem1} and Proposition~\ref{prop:appl2-lem2}.
\end{proof}

As we did for the alternative Daugavet property in
Corollary~\ref{cor:ADP-SCD-operators-nonseparable}, we can extend
Theorem~\ref{appl2-the1} to the non separable case.

\begin{corollary}
Let $X$ be a Banach space with the Daugavet property and let $T\in
L(X)$ be such that $T|_Y$ is an hereditary-SCD-operator for every
separable subspace $Y$ of $X$. Then, $T$ is narrow and, in
particular, $\|\Id+ T\|=1+\|T\|$.
\end{corollary}

\begin{proof}
We fix $x,y\in S_X$, a slice $S$ of $B_X$ and $\eps>0$. We take a
separable subspace $Y_1$ of $X$ such that $x,y\in Y_1$ and such
that $S\cap Y_1\neq \emptyset$ and we follow the proof of
Corollary~\ref{cor:ADP-SCD-operators-nonseparable}, using
\cite[Theorem~4.5]{KadSW2} instead of
Remark~\ref{rem:ADP-separably-determined}, to get a separable
subspace $Y$ of $X$, $T$-invariant, with the Daugavet property and
such that $x,y\in S_X$ and $S\cap Y\neq \emptyset$. Now, $T|_Y$ is
a hereditary-SCD-operator, so Theorem~\ref{appl2-the1} gives us
that $T|_Y$ is narrow. Then, we may find $z\in S\cap Y\subseteq S$
such that $\|x+z\|\geq 2-\eps$ and
\begin{equation*}
\|Ty-Tz\|=\|T_Y(y) - T|_Y(z)\| < \eps.\qedhere
\end{equation*}
\end{proof}

The following particular cases are especially interesting. The
first one was proved in \cite[Theorem~4.13]{KadSW2} with a
different argument.

\begin{corollary}
Let $X$ be a Banach space with the Daugavet property and let $T\in
L(X)$ be an operator which does not fix copies of $\ell_1$. Then,
$T$ is narrow.
\end{corollary}

\begin{corollary}
Let $X$ be a Banach space with the Daugavet property and let $T\in
L(X)$ be an operator such that $T(B_X)$ is strongly regular. Then,
$T$ is narrow.
\end{corollary}

\begin{remarks}$ $\label{remark-rightideal}
\begin{enumerate}
\item[(a)] {\slshape The class of hereditary-SCD-operators is
    a right operator ideal.\ } Indeed, if
    $T:X_1\longrightarrow X_2$ is an arbitrary operator and
    $S:X_2\longrightarrow X_3$ is a hereditary-SCD-operator,
    then $[ST](B_{X_1})\subseteq S(\|T\|B_{X_2})$, so $ST$ is
    an hereditary-SCD-operator.
\item[(b)] {\slshape The class of hereditary-SCD-operators is
    not a left operator ideal.\ } Indeed, we consider a
    norm-one projection $T:L_1[0,1]\longrightarrow X\equiv
    \ell_1$ which is a hereditary-SCD-operator since $\ell_1$
    is RNP. We also consider a quotient map
    $S:\ell_1\longrightarrow \ell_1/Y\equiv L_1[0,1]$ (by just
    using the factor universality of $\ell_1$). Then,
    $ST(B_{L_1[0,1]})=B_{L_1[0,1]}$ so $ST$ is not even an
    SCD-operator.
\item[(c)] As a consequence, {\slshape there are narrow
    operators which are not SCD-operators.\ } Indeed, since
    the set of narrow operators is clearly a left operator
    ideal, the operator $ST$ above is narrow.
\end{enumerate}
\end{remarks}

\section{Countable $\pi$-bases of the weak topology}\label{sec:countable-pi-base}
It was shown in Proposition~\ref{SCD:prop-countablepibase} that a
convex bounded subset $A$ of a Banach space $X$ is SCD if it has a
countable $\pi$-base of the weak topology. But we do not know
whether these two properties are equivalent. The aim of this
section is to discuss this possible equivalence. In a first
subsection we will show that the class of sets having countable
$\pi$-bases of the weak topology contains separable CPCP sets. We
already know that it contains those sets which do not have
$\ell_1$-sequences (Theorem~\ref{Examples-Theorem:ell1}), so this
class covers most of the examples of SCD sets presented in this
paper. In the second subsection we will show that convex bounded
subsets of both $\ell_1(c_0)$ and $c_0(\ell_1)$ also have
countable $\pi$-bases of the weak topology. Finally, the third
subsection contains several characterizations of SCD sets which
remind of the property we are dealing with.

\subsection{CPCP sets}
We start with a sufficient condition to have a countable
$\pi$-base of the weak topology.

\begin{prop}\label{prop:dense=>piseparable}
Let $X$ be a Banach space and let $A$ be a separable closed convex
bounded subset of $X$ such that there is a weakly dense subset $B$
of $A$ consisting of points of continuity of
$\,\Id:(A,\sigma(X,X^*))\longrightarrow (A,\|\cdot\|)$. Then,
$(A,\sigma(X,X^*))$ has a countable $\pi$-base.
\end{prop}

\begin{proof}
Let $D$ be a countable norm dense subset of $B$, and for every
$d\in D$ and every $n\in\mathbb{N}$ let $U_d^n$ be a weak open
neighborhood of $d$ in $A$ of diameter less than $\frac{1}{n}$. We
claim that the countable family $\{U_d^n\,:\,n\in\N,\,d\in D\}$ is
a $\pi$-base of $A$. Indeed, let $W$ be a weakly open subset of
$A$. Since $B$ is weakly dense in $A$, $W\cap B$ is non-empty and
relatively norm open in $B$ so, since $D$ is norm dense in $B$,
there is $d\in D\cap W$. Now, $W$ is a norm open neighborhood of
$d$ relative to $A$, so it contains $B(d,1/n)\cap A$ for some
$n\in \mathbb{N}$ and so $U_d^n\subseteq W$. We are done.
\end{proof}

A first consequence of the above result deals with LUR renorming.
It is clear from the definition that denting points are points of
weak-norm continuity of the identity map and so, as it was
commented before Example~\ref{example:LUR=>SCD}, the unit ball of
a Banach space with a LUR norm fulfills the above condition. It
was also commented there that every separable Banach space can be
equivalently renormed with a LUR norm.

\begin{example}$ $
\begin{enumerate}
\item[(a)] {\slshape Let $X$ be a separable Banach space with
    a LUR norm. Then, $B_X$ has a countable $\pi$-base of the
    weak topology.\ }
\item[(b)] As a consequence, {\slshape every separable Banach
    space $X$ admits an equivalent norm $|\cdot|$ such that
    $B_{(X,|\cdot|)}$ has a countable $\pi$-base of the weak
    topology.\ }
\end{enumerate}
\end{example}

We are going to show that CPCP sets have countable $\pi$-bases for
the weak topology. We recall that a closed convex bounded subset
$A$ of a Banach space $X$ has the CPCP if every convex closed
subset $B$ of $A$ contains a weak-to-norm point of continuity of
the identity mapping. In this case, for every convex subset $B$ of
$A$ and for every $\eps>0$, there is a relatively weakly open
subset $C\subseteq B$ with $\diam(C)<\eps$ \cite{Bourgain1980}. We
need the following result which follows from
\cite[Lemma~I.0]{Gho-God-Mau-Scha}; we haven't found a direct
reference, so we include a proof for the sake of completeness.

\begin{lemma}
Let $X$ be a Banach space and let $A$ be a closed convex bounded
subset of $X$ with the CPCP. Then, there is a weakly dense subset
$D$ of $A$ consisting of points of weak-norm continuity of
$\Id:(A,\sigma(X,X^*))\longrightarrow (A,\|\cdot\|)$.
\end{lemma}

\begin{proof}
We fix a sequence of positive $\eps_n$ tending to zero and write
$$
D_n=\bigcup \{C\ : \ C\ \text{is weakly open in $A$ and $\diam(C)<\eps_n$}\}.
$$
Let us prove that $D=\bigcap_{n\in\N} D_n$ is weakly dense in $A$.
Indeed, let $U\subseteq A$ be relatively weakly open. We pick
$U_1\subseteq U$ convex closed with non-empty interior. Then,
there is a relatively weakly open subset $C_1$ of $A$ of diameter
less than $\eps_1$ such that $C_1$ is contained in the weak
interior of $U_1$. We repeat the process to find a decreasing
sequence $C_n$ of weakly open sets with non-empty interior such
that $\diam(C_n)<\eps_n$ and $\overline{C_{n+1}}\subseteq C_{n}$.
Then, the Cantor theorem tells us that there is $x\in
\bigcap_{n\in \N} \overline{C_n}$. Now, we have in particular that
$x\in C_1\subseteq U_1\subseteq U$. On the other hand, for every
$n\in\N$, $x\in C_n$ and $\diam(C_n)<\eps_n$, so $x\in D_n$.
Therefore, $x\in D$. Finally, every point of $D$ has weak
neighborhoods of arbitrarily small diameter, showing that it is a
point of continuity.
\end{proof}

This result, together with
Proposition~\ref{prop:dense=>piseparable} gives the main result of
the subsection.

\begin{corollary}
Let $X$ be a Banach space and let $A$ be a separable closed convex
bounded subset of $X$ with the CPCP. Then, $A$ has a countable
$\pi$-base for the weak topology.
\end{corollary}

With the above result, most of the types of SCD sets presented in
the section~\ref{sec:sets} have a countable $\pi$-base of the weak
topology. The only exception is the family of strongly regular
sets which are not CPCP. There are two main examples of sets of
this kind, but in both cases, the sets have a countable $\pi$-base
of the weak topology.

\begin{examples}$ $
\begin{enumerate}
\item[(a)] {\slshape The set constructed by S.~Argyros,
    E.~Odell, and H.~Rosenthal \cite{ArgOdeRos} which is
    strongly regular but does not have the CPCP is a subset of
    $c_0$, so it has a countable $\pi$-base of the weak
    topology since it does not have $\ell_1$-sequences.\ }
\item[(b)] {\slshape The set constructed by W.~Schachermayer
    \cite{Scha-example} which is a subset $C$ of a Banach
    space $Z$ which does not have the CPCP but $Z^{**}$ is
    strongly regular (so $Z$ is strongly regular). But then,
    $(C,\sigma(X,X^*))$ has a countable $\pi$-base of the weak
    topology since $Z$ does not contain $\ell_1$.\ }
\end{enumerate}
\end{examples}

\subsection{$\boldsymbol{c_0(\ell_1)}$ and $\boldsymbol{\ell_1(c_0)}$} Our goal in this subsection is to show that convex bounded subsets
of the spaces $c_0(\ell_1)$ and $\ell_1(c_0)$ have a countable
$\pi$-base of the weak topology. The first case is easier to
demonstrate.

\begin{example}
{\slshape Every bounded convex subset $A$ of the space
$c_0(\ell_1)$ has a countable $\pi$-base of the weak topology.\ }
\end{example}

\begin{proof}
Let $X$ denote $c_0(\ell_1)$. For every $m\in \N$, we denote
$$
Y_m=\bigl[\ell_1\oplus\ell_1\oplus \overset{m}{\ldots}\oplus \ell_1
\oplus 0 \oplus 0 \oplus \cdots\bigr]_\infty \subseteq c_0(\ell_1)
$$
and $P_m:X\longrightarrow Y_m$ for the natural projection. Since
$P_m(A)$ is a convex bounded subset of $Y_m$ and $Y_m$ is
isomorphic to $\ell_1$, there is a countable $\pi$-base
$\{{S}_{m,k}\,:\, k \in \N\}$ of $(P_m(A), \sigma(Y_m,Y_m^*))$. We
are going to prove that the collection
$$
\widetilde{S}_{m,k}=
\bigl[P_m^{-1}({S}_{m,k})\bigr] \cap A \qquad (m,k \in \N)
$$
forms a countable $\pi$-base of $(A, \sigma(X,X^*))$. Indeed, let
$U, V$ be weak neighborhoods of $0$, $V + V \subseteq U$, $a \in
A$, and denote $B = (a + U) \cap A$. Every relatively weakly open
subset of $A$ is of the same form as $B$, so we have to prove that
$\widetilde{S}_{m,k} \subseteq B$ for some choice of $m$ and $k$.
Assume to the contrary that none of $\widetilde{S}_{m,k}$ is
contained in $B$. For $m\in\N$ big enough, all the $P_m(A)$
intersect $(a + V)$. Fix $m \in \N$ with $C_m = (a + V) \cap
P_m(A) \neq \emptyset$. Then there is $k(m) \in \N$ with
${S}_{m,k(m)} \subseteq C_m$. According to our assumption
$\widetilde{S}_{m,k(m)}$ is not contained in $B$, so there is an
$x_m \in \widetilde{S}_{m,k(m)} \setminus B$. This $x_m$ can be
written as $x_m = y_m + z_m$, where $y_m \in {S}_{m,k(m)}
\subseteq C_m$ and $z_m \in \Ker P_m$. Since $x_m \in A$ and $y_m
\in P_m(A)$, we have that $z_m$ is a bounded sequence, and since
by our construction $(z_m)$ tend to $0$ coordinate-wise as $m \to
\infty$, we can deduce that $(z_m) \longrightarrow 0$ in the weak
topology. Therefore, for some $m$ big enough $z_m \in V$ and
consequently $x_m = y_m + z_m \in (a + V) + V \subseteq a + U$.
Since $x_m \in A$, this means that $x_m \in (a + U) \cap A = B$,
which contradicts the selection of $x_m$.
\end{proof}

\begin{remark}
{\slshape The argument above also works for $c_0$-sums of RNP
spaces.\ } Indeed, this follows from the fact that a finite-sum of
RNP spaces is again a RNP space (see \cite[Theorem~6.5.b]{CasGon},
for instance).
\end{remark}

Let us remark with an example that to have a countable $\pi$-base
of the weak topology does not imply that any point has a countable
base of weak neighborhoods.

\begin{example}
{\slshape The unit ball of $X=c_0(\ell_1)$ has no point with a
countable base of relative weak neighborhoods.\ } Indeed, we
consider an arbitrary $x = (x_n)_{n \in \N}\in B_X$, where $x_n
\in \ell_1$, $\|x_n\| \longrightarrow 0$ and $\max_{n \in
\N}\|x_n\| \leq 1$. We fix $n_0\in\N$ such that $\|x_{n_0}\| <
1/2$ and we consider the subset
$$
A=\left\{(y_n)_{n\in\N}\in B_X\,:\ y_n=x_n \text{ if $n\neq n_0$},\
\|x_{n_0}-y_{n_0}\|\leq 1/2\right\}.
$$
Then, $A$ is a closed subset of $B_X$ containing $x$, so if $x$
has a countable base of relative weak neighborhoods in $B_X$, then
$x$ has also a countable base of relative weak neighborhoods in
$A$. But the latter is impossible, because $A$ is affinely
homeomorphic to $B_{\ell_1}$, with $x$ being the image of $0 \in
B_{\ell_1}$.
\end{example}

To get the second example we need a technical result.

\begin{lemma}\label{lemma:local-pi-base-local-eps-base}
Let $X$ be a separable Banach space. Then, the following are
equivalent.
\begin{enumerate}
\item[(i)] Every convex bounded subset of $X$ has a countable
    $\pi$-base of the weak topology.
\item[(ii)] Every closed convex bounded subset $A$ of $X$ has
    a point with a countable local $\pi$-base of relatively
    weakly open subsets (i.e.\ there is $x\in A$ and a
    sequence $\{U_n\,:\,n\in\N\}$ of relatively weakly open
    subsets of $A$ such that for every relative weak
    neighborhood $V$ of $x$ there is some $U_n \subseteq V$.)
\item[(iii)] For every $\eps > 0$, every closed convex bounded
    subset $A$ of $X$ has a point with a countable local
    $\eps$-base of relatively weakly open subsets (i.e.\ there
    is $x\in A$ and a sequence $\{U_n\,:\,n\in\N\}$ of
    relatively weakly open subsets of $A$ such that for every
    weakly open neighborhood $V \subseteq X$ of $x$ \emph{in
    the whole space} there is $n \in \N$ with $U_n \subseteq V
    + \eps B_X$.)
\end{enumerate}
\end{lemma}

\begin{proof}
(i) $\Rightarrow$ (ii) $\Rightarrow$ (iii) are clear since a
$\pi$-base is a local $\pi$-base, and a local $\pi$-base is an
$\eps$-base for every $\eps>0$.

(iii) $\Rightarrow$ (i). It is straightforward to show that it is
enough to deal with \emph{closed} convex bounded subsets of $X$.
Just observe that if $\{U_n\,:\,n\in\N\}$ is a $\pi$-base for the
weak topology of the closure of a bounded convex subset $A$ of
$X$, then $\{U_n\cap A\,:\,n\in\N\}$ is a $\pi$-base of the weak
topology of $A$ itself.

We then fix a \emph{closed} convex bounded subset $A \subseteq X$.
We first remark that, for every $\eps>0$, the subset $B_\eps
\subseteq A$ of points having a countable local $\eps$-base is
weakly dense in $A$. Indeed, we consider an arbitrary weakly open
subset $U$ of $X$ intersecting $A$ and we fix another weakly open
subset $V \subseteq U \subseteq X$ intersecting $A$ with
$\overline{V}^{\sigma(X,X^*)} \subseteq U$. According to our
assumption, there is $x \in \overline{V}^{\sigma(X,X^*)}\cap A$
with a local $\eps$-base $\{U_n\,:\,n\in\N\}$. But then, $V_n= U_n
\cap V$ form a countable local $\eps$-base of relatively weakly
open subsets of $A$ for $x$, i.e.\ $x \in B_\eps \cap U$, so $B
\cap U \neq \emptyset$.

Now, for every $k\in\N$ we take a countable norm dense subset
$\{b_{k,m}\,:\,m \in \N\}$ in $B_{1/k}$, and for every $b_{k,m}$,
we select a $1/k$-base $\{U_{k,m,n}\,:\,n \in \N\}$. Let us show
that $\{U_{k,m,n}\,:\,k,n,m \in \N\}$ forms a $\pi$-base for $(A,
\sigma(X,X^*))$. Indeed, let $U, V$ be weak neighborhoods of $0$,
$V + V \subseteq U$, $a \in A$, and denote $G = (a + U) \cap A$.
We have to prove that $U_{k,m,n} \subseteq G$ for some choice of
$k,m,n\in\N$. To do this, we take $k\in \N$ big enough that
$\frac1k B_X \subseteq V$. According to our construction, there is
$m\in\N$ with $b_{k,m} \in (a + V) \cap A$. Then, there is
$n\in\N$ with $U_{k,m,n} \subseteq (a + V) + \frac1k B_X$.
Therefore
\begin{equation*}
U_{k,m,n} \subseteq \left(a + V + \frac1k B_X\right) \cap A
\subseteq (a + V + V) \cap A \subseteq (a + U) \cap A = G.\qedhere
\end{equation*}
\end{proof}

We are now able to present the second example.

\begin{example}
{\slshape Every bounded convex subset $A$ of the space
$\ell_1(c_0)$ has a countable $\pi$-base of the weak topology.\ }
\end{example}

\begin{proof}
Let $X$ denote $\ell_1(c_0)$. For $\eps>0$ fixed, arguing the same
way as in the beginning of the proof of
Theorem~\ref{Examples-th4}, we select an open slice $S \subseteq
A$ and an $m\in \N$ with the following property
\begin{equation}\label{eq:PROPERTY-c_0-ell_1}
\bigl\|(0, \ldots, 0, y_{m+1}, y_{m+2}, \ldots)\bigr\|
<\frac{\eps}{2} \qquad \bigl((y_n)_{n\in \N} \in S\bigr).
\end{equation}
Let us prove that every $x_0 \in S$ has a countable $\eps$-base of
relatively weakly open subsets and
Lemma~\ref{lemma:local-pi-base-local-eps-base} will give the
result.

We denote
$$
Y_m=\bigl[c_0\oplus c_0\oplus \overset{m}{\ldots}\oplus c_0 \oplus 0
\oplus 0 \oplus \cdots\bigr]_{\ell_1} \subseteq \ell_1(c_0)
$$
and let $P_m:X\longrightarrow Y_m$ be the natural projection.
Since $Y_m$ is isomorphic to $c_0$, there is a countable local
$\pi$-base $\{U_n\,:\,n \in \N\}$ of $P_m(x_0)$ in $(P_m(S),
\sigma(X,X^*))$. Consider
$$
\widetilde{U_n} = P_m^{-1}(U_n) \cap S \qquad (n\in\N)
$$
which are weakly open subsets of $S$, and hence they are weakly
open in $A$. Let us show that $\{\widetilde{U_n}\,:\,n\in \N\}$
forms an $\eps$-base for $x_0$ in $A$. Consider a weakly open
neighborhood $V \subseteq X$ of $x_0$. By
\eqref{eq:PROPERTY-c_0-ell_1}, we have that
\begin{equation} \label{eq1ps}
\|P_m(y) - y\| < \eps /2 \qquad (y \in S).
\end{equation}
So $\left(V + \frac\eps 2 B_X\right) \cap P_m(S)$ is a weak
neighborhood of $P_m(x_0)$ in $P_m(S)$. So there is an $n\in\N$
such that $U_n \subseteq V + \frac\eps 2 B_X$. Applying
\eqref{eq1ps} once more, we obtain that
\begin{equation*}
\widetilde{U_n} = P_m^{-1}(U_n) \cap S \subseteq U_n + \frac\eps2
B_X \subseteq V + \eps B_X. \qedhere
\end{equation*}
\end{proof}

\subsection{Two characterizations of SCD sets}
The aim of this part of the section is to establish some
characterizations of SCD sets which remind of countable
$\pi$-bases of the weak topology. The first one deals with convex
combinations of slices.

\begin{theorem} \label{weakly regular}
A bounded convex subset $A$ of a Banach space $X$ is an SCD set if
and only if there is a sequence $\{V_n\,:\,n\in\N\}$ of convex
combinations of slices of $A$ such that every relatively weakly
open subset of $A$ contains some of the $V_n$.
\end{theorem}

\begin{proof}
The ``if'' part is direct consequence of Propositions
\ref{determ-sequences} and \ref{SCD-sets-rem-2}.

Conversely, asume that $A$ is an SCD set and suppose without loss
of generality that $A \subseteq B_X$. Let $S_n = S(A,x_n^*,
\eps_n)$ for $n \in \N$, be a determining sequence of slices for
$A$. Let us show that the convex combinations of the $S_n$'s with
rational coefficients form the countable collection of convex
combinations of slices that we need. Indeed, let $U$ be a
relatively weakly open subset of $A$. Select another relatively
weakly open subset $V \subseteq U$ such that $\alpha=\dist(V, A
\setminus U) >0$. Due to Bourgain's lemma
(Lemma~\ref{lemma:Bourgain}), there is a convex combination of
slices $\sum_{j=1}^{m}\lambda_j G_j\subseteq V$. According to
Proposition~\ref{determ-sequences}, for every $j = 1,2, \ldots ,
m$ there is $n(j)\in\N$ such that $S_{n(j)} \subseteq G_j$. Then,
$\sum_{j=1}^{m}\lambda_j S_{n(j)} \subseteq V$. What remains is to
find rationals $\mu_j > 0$ with $\sum_{j=1}^{m}\mu_j = 1$ and
$|\mu_j - \lambda_j| < \alpha$. Then, the Hausdorff distance
between $\sum_{j=1}^{m}\mu_j S_{n(j)}$ and
$\sum_{j=1}^{m}\lambda_j S_{n(j)}$ is smaller than $\alpha$, so
$\sum_{j=1}^{m}\mu_j S_{n(j)} \subseteq V + \alpha\,B_X\subseteq
U$.
\end{proof}

The second result gives a reformulation of SCD in terms of
topological properties of the set of extreme points of its
weak$^*$ closure in the bidual. For a convex bounded subset $A$ of
a Banach space $X$, denote $\overline{A}^{**}$ the weak-star
closure of $A$ in $X^{**}$.

\begin{theorem} \label{SCD sets-th2}
Let $X$ be a Banach space and let $A$ be a convex bounded subset
of $X$. Put $W=\left(\extr\bigl(\overline{A}^{**}\bigr),
\sigma(X^{**}, X^{*}) \right)$. Then, the following are
equivalent:
\begin{enumerate}
\item[(i)] $A$ is an SCD set.
\item[(ii)] $W$ has a countable $\pi$-base.
\end{enumerate}
\end{theorem}

\begin{proof}
(i) $\Longrightarrow$ (ii). We take a sequence of slices $S_n =
S(A,x_n^*, \eps_n)$ for $n \in \N$ which is determining for $A$
and we write
$$
S_n^{**} = S\bigl(\overline{A}^{**},x_n^*, \eps_n\bigr) \subseteq
\overline{A}^{**}
$$
for the natural extensions of $S_n$ to slices of
$\overline{A}^{**}$. Then, the family $U_n= S_n^{**} \cap W$ for
$n \in \N$ forms a $\pi$-base of $W$. Indeed, we consider a
relatively weak$^*$-open subset $U$ of $W$. Due to Choquet's lemma
(that for any locally convex topology, slices containing an
extreme point of a compact convex set make up a neighborhood base
of the extreme point, see \cite[Definition~25.3 and
Proposition~25.13]{ChoLecAnal-II}), there is a slice $S^{**}=
S\bigl(\overline{A}^{**},x^*, \eps\bigr)$ of $\overline{A}^{**}$
generated by some $x^* \in X^*$ and $\eps>0$ such that $U
\supseteq S^{**} \cap W \neq \emptyset$. Now, according to
Proposition~\ref{determ-sequences}, there is an $n \in \N$ such
that
$$
S_n \subseteq S(A,x^*, \eps/2) \subseteq
S\bigl(\overline{A}^{**},x^*, \eps/2\bigr).
$$
Then, $S_n^{**}$ is contained in the relative weak$^*$-closure of
$S\bigl(\overline{A}^{**}, x^*, \eps/2\bigr)$ in
$\overline{A}^{**}$, so $S_n^{**} \subseteq S^{**}$ and
 $$
 U_n= S_n^{**} \cap W \subseteq S^{**}\cap W \subseteq U.
 $$

(ii) $\Longrightarrow$ (i). We consider a countable $\pi$-base
$\{U_n\,:\, n\in \N\}$ of $W$ consisting of relatively
weak$^*$-star open subsets. Again by Choquet's lemma, there are
$x_n^* \in X^*$ and $\eps_n > 0$ such that
$$
U_n \supseteq \widetilde{U_n}= S\bigl(\overline{A}^{**},x_n^*,
\eps_n\bigr) \cap W \neq \emptyset.
$$
Let us prove that the slices $S_{n,m} = S(A,x_n^*, 1/m)$ with $n,m
\in \N$, form a determining sequence for $A$. Indeed, we denote
$S_{n,m}^{**}$ the \emph{closed} slices of $\overline{A}^{**}$
generated by $x_n^*$ and $1/m$. For every slice $S=S(A,x^*, \eps)$
of $A$, since $\{\widetilde{U_n}\,:\,n\in \N\}$ is a $\pi$-base of
$W$, there is $n \in \N$ such that
$$
S^{**} \cap W  \supseteq \widetilde{U_n} \qquad \text{where }
S^{**}=S\bigl(\overline{A}^{**},x^*, \eps_n\bigr),
$$
so for $m\in \N$ big enough we have
$$
S^{**} \cap W \supseteq S_{n,m}^{**} \cap W.
$$
Then, taking into account that, for every $n\in\N$,
$$
G_n= \bigcap_{m\in \N} S_{n,m}^{**}
$$
is a closed face of $\overline{A}^{**}$, the Krein-Milman theorem
gives us that
$$
G_n=\overline{\conv\bigl(G_n \cap
W\bigr)}^{\sigma(X^{**},X^*)}.
$$
Therefore,
$$
S^{**} \supseteq \overline{\conv\bigl(S^{**} \cap
W\bigl)}^{\sigma(X^{**},X^*)} \supseteq
\overline{\conv\left(\bigcap_{m\in \N} S_{n,m}^{**} \cap
W\right)}^{\sigma(X^{**},X^*)} = G_n.
$$
This means that the intersection of the decreasing sequence of
$\sigma(X^{**},X^*)$ compact sets $\{S_{n,m}^{**}\,:\,m\in\N\}$ is
contained in $S^{**}$. But $S^{**}$ is a relatively
$\sigma(X^{**},X^*)$ open set in $\overline{A}^{**}$, so for
sufficiently big $m\in\N$, all the $S_{n,m}^{**}$ are subsets of
$S^{**}$. For these $m$, we have
$$
S = S^{**} \cap A \supseteq S_{n,m}^{**}  \cap A \supseteq S_{n,m}.
$$
Finally, we use the characterization of SCD sets from
Proposition~\ref{determ-sequences}.
\end{proof}

The following is an easy consequence of the above result.

\begin{corollary}
Let $X$ be a Banach space and let $A$ be a bounded convex subset
of $X$. If $A$ is SCD, then
$\left(\extr\bigl(\overline{A}^{**}\bigr),\sigma(X^{**},X^*)\right)$
is separable.
\end{corollary}

The particular case of the above corollary for subsets of
separable Banach spaces without copies of $\ell_1$ should be
previously known. Anyway, we include an easy direct proof of this
fact.

\begin{remark}
{\slshape Let $X$ be a separable Banach space without copies of
$\ell_1$ and let $A$ be a convex bounded subset of $X$. Then,
$\left(\extr\bigl(\overline{A}^{**}\bigr),\sigma(X^{**},X^*)\right)$
is separable.\ } Indeed, we write
$$
C=\conv\bigl(\extr\bigl(\overline{A}^{**}\bigr)\bigr)
$$
and we observe that $C$ is $\sigma(X^{**},X^*)$-sequentially dense
in its weak$^*$-closure $\overline{A}^{**}$ (see
\cite[Theorem~4.1]{van-Duslt}). Then, we take a sequence
$\{y_n\,:\,n\in\N\}$ dense in $A$ and we consider those extreme
points of $\overline{A}^{**}$ needed to approximate each $y_n$ by
a sequence of convex combinations. The union of all these extreme
points (while countable) is weak$^*$-dense in the set of all
extreme points of $\overline{A}^{**}$ by the reversed Krein-Milman
theorem.
\end{remark}

\section{Open questions}\label{sec:openquestions}

\begin{question}\label{question-1}
Let $X$ be a Banach space and let $A$ be a convex bounded subset
of $X$. If $A$ is SCD, does $A$ have a countable $\pi$-base for
the weak topology?
\end{question}

\begin{question}\label{question-2}
Let $X$ be an SCD space. Does every convex bounded subset of $X$
have a countable $\pi$-base for the weak topology?
\end{question}

Related to these questions is the following one.

\begin{question}
Let $L$ be a compact subset of a locally convex space and let $K$
be its closed convex hull. If $L$ has a countable $\pi$-base, does
it imply that $K$ also has a countable $\pi$-base? What if
$L=\extr(K)$?
\end{question}

Let us explain why this question is related to the above two.
Observe that if $D$ is a dense subspace of a topological space $E$
and $\mathcal{B}$ is a $\pi$-base for $E$, then $\{B\cap D \,:\,
B\in\mathcal{B}\}$ is a $\pi$-base for $D$. In particular, if
$(\overline{A}^{\ast\ast},\sigma(X^{\ast\ast},X^\ast))$ has a
countable $\pi$-base, then so does $(A,\sigma(X,X^\ast))$. Thus, a
positive answer to the preceding question combined with
Theorem~\ref{SCD sets-th2} would imply a positive answer to
Questions \ref{question-1} and \ref{question-2}.

\begin{questions}$ $
\begin{enumerate}
\item[(a)] Is every Banach space with unconditional basis SCD?
\item[(b)] A simpler case: let $X$ be a Banach space with
    $1$-symmetric basis. Is $B_X$ an SCD set?
\end{enumerate}
\end{questions}

\begin{question}
Are the concepts of SCD sets and almost-SCD sets equivalent (see
Remark~\ref{remark:ADP=lush_with_almostSCD} for the definition)?
\end{question}

\begin{questions}
Let $X$ be a separable Banach space such that no subspace of it
can be renormed with the Daugavet property. Is $X$ SCD?
\end{questions}

\begin{questions}$ $
\begin{enumerate}
\item[(a)] Is the sum of two SCD-operators an SCD-operator?
\item[(b)] Is the sum of two hereditary-SCD-operators a
    hereditary-SCD-operator?
\end{enumerate}
\end{questions}

\noindent \textbf{Acknowledgments.\ } We would like to thank Bill
Johnson and Rafael Pay\'{a} for useful comments on the subject of the
paper. We also thank the referee for the carefully reading of the
manuscript and for multiple stylistic remarks.

\end{document}